\documentclass[12pt]{amsart}
\usepackage[margin=1in]{geometry}
\usepackage{amsaddr}
\usepackage[toc,page]{appendix}
\usepackage{graphicx}
\usepackage{epsfig, epstopdf}
\usepackage{multicol}
\usepackage{psfrag}
\usepackage{amssymb}
\usepackage{amsthm}
\usepackage{lineno}
\usepackage{amsfonts}
\usepackage{amsmath}
\usepackage{caption}
\usepackage{subcaption}
\usepackage{bm}
\usepackage{booktabs}  
\usepackage{enumitem}
\usepackage{mathtools}
\usepackage{flafter}

\usepackage{framed}
\usepackage{color}
\usepackage{hyperref}
\usepackage[english,algosection]{algorithm2e} 
\hypersetup{colorlinks,breaklinks,
           linkcolor=blue,urlcolor=blue,
           anchorcolor=blue,citecolor=blue}
\usepackage[normalem]{ulem}
\usepackage{float}
\usepackage{comment}
\usepackage{setspace}
\newtheorem{thm}{Theorem}[section]

\newtheorem{property}{Property}[section]
\newtheorem{lemma}{Lemma}[section]

\newtheorem{remark}{Remark}
\graphicspath{{fig/}}

\newcommand{\V}{\mathbb{V}}

\newcommand{\TT}{\mathcal{T}}

\newcommand{\R}{\mathbb{R}}
\newcommand{\N}{\mathcal{N}}
\newcommand{\E}{\mathbb{E}}

\newcommand{\bn}{\mathbf{n}}
\newcommand{\bq}{\mathbf{q}}

\newcommand{\f}{\mathbf{f}}

\newcommand{\s}{\mathbf{s}}

\newcommand{\x}{\mathbf{x}}

\newcommand{\y}{\mathbf{y}}

\newcommand{\Q}{\mathcal{Q}}

\newcommand{\brho}{\bm \rho}

\newcommand{\norm}[1]{\left\Vert#1\right\Vert}

\newcommand{\bxi}{ \bm{\xi} }
\newcommand{\slk}{k}

\def\bR{{\mathbf R}}
\newcommand{\T}{{\mathcal T}}
\def\bzeta{{\boldsymbol \zeta}}
\def\zetacoeff{{\zeta}}
\def\bW{{W}}

\def\bPi{{\Pi}}

\DeclarePairedDelimiter{\ceil}{\lceil}{\rceil}

\newcommand{\bb}{b}
\newcommand{\bu}{u}

\renewcommand{\div}{\operatorname{div}}
\hypersetup{
colorlinks=true,
linkcolor=blue,
citecolor=red}
\definecolor{ogreen}{rgb}{0.2, 0.8, 0.6}
\definecolor{dgreen}{rgb}{0.3098,0.3843,0.1569}
\definecolor{blue1}{RGB}{176,196,222 }
\definecolor{blue2}{RGB}{54,100,139}
\definecolor{blue3}{RGB}{80,120,170}

\newcommand{\hnote}[1]{}
\newcommand{\hnotehelp}[1]{}

\newcommand{\urev}[1]{\textcolor{black}{#1}}

\newcommand{\hrev}[1]{\textcolor{black}{#1}}
\newcommand{\rsout}[1]{}

\title[ML Hierarchical Decomposition of Finite Element White Noise]{Multilevel Hierarchical Decomposition of Finite Element White Noise with Application to Multilevel Markov Chain Monte Carlo}
\author[]{Hillary R. Fairbanks$^{1}$ \and Umberto Villa$^{2}$\and Panayot S. Vassilevski$^{1,3}$}
\address{$^{1}$Center for Applied Scientific Computing, Lawrence Livermore National Laboratory, Livermore, CA, USA \\ $^{2}$Electrical \& Systems Engineering Department, Washington University in St. Louis, St. Louis, MO, USA\\ $^{3}$Fariborz Maseeh Department of Mathematics and Statistics, Portland State University, Portland, OR, USA}

 \thanks{This work was performed under the auspices of the U.S. Department     
 of Energy by Lawrence Livermore National Laboratory under Contract             
 DE-AC52-07NA27344 (LLNL-JRNL-820098)}

\begin{document}
\maketitle

\begin{abstract}
In this work we develop a new hierarchical multilevel approach to generate Gaussian random field realizations in \hrev{an algorithmically} scalable manner that is well-suited to incorporate into multilevel Markov chain Monte Carlo (MCMC) algorithms. This approach builds off of other partial differential equation (PDE) approaches for generating Gaussian random field realizations; in particular, a single field realization may be formed by solving a reaction-diffusion PDE with a spatial white noise source function as the righthand side. While these approaches have been explored to accelerate forward uncertainty quantification tasks, e.g. multilevel Monte Carlo, the previous constructions are not directly applicable to multilevel MCMC frameworks which build fine scale random fields in a hierarchical fashion from coarse scale random fields. Our new hierarchical multilevel method relies on a hierarchical decomposition of the white noise source function in $L^2$ which 
allows us to form Gaussian random field realizations across multiple levels of discretization in a way that fits into multilevel MCMC algorithmic frameworks. 
After presenting our main theoretical results and numerical scaling results to showcase the utility of this new hierarchical PDE method for generating Gaussian random field realizations, this method is tested on a \hrev{four}-level MCMC algorithm to explore its feasibility. 
\end{abstract}
\section{Introduction}
Spatially correlated random fields are commonly used in the numerical simulation of partial differential equations (PDEs) with variable coefficients. In the case where these coefficients are not well known,
as is typically the case in many geophysics applications where the coefficient describes a physical parameter, the coefficient is modeled as a random field,
and uncertainty quantification (UQ) may be applied as a tool to assess the reliability of the model as well as the sensitivity to changes in this parameter. To reduce the uncertainty in the system, we may further improve the model by utilizing observational data in a Bayesian framework. That is, data related to the model output, as well as information about the model may be combined to learn the probability distribution of the variable coefficient.

For large-scale applications, common methods to perform Bayesian inference are infeasible. With the refinement of the spatial discretization scheme, both forming realizations of these random fields and performing forward PDE simulations are computationally demanding, as many approaches do not scale with the increase in problem size. Furthermore, Bayesian inference approaches are typically limited to Markov chain Monte Carlo (MCMC)~\cite{metropolis1953equation,hastings1970monte,robert2013monte}, and its variants, which require a large number of simulations as the parameter space is explored. However infeasible this approach may be, MCMC methods still lie at the root of many nonlinear Bayesian inference algorithms due to ease of implementation as well as the ability to be applied in a blackbox fashion.

Over the recent decades new MCMC approaches have been developed to accelerate the parameter space exploration, in some cases allowing to perform nonlinear Bayesian inference on large-scale applications. Notable approaches include those which utilize local approximations of the Hessian and gradient to modify the MCMC proposal~\cite{martin2012stochastic,petra2014computational,cui2016dimension,CuiMartinMarzoukEtAl14,Bui-ThanhGirolami14,BeskosGirolamiLanEtAl17}. This method, dubbed, stochastic Newton MCMC, was shown in \cite{petra2014computational} to accelerate mixing; however, it requires gradient and Hessian information in addition to having solvers for the forward PDE and adjoint PDE models, as opposed to simply having the forward PDE model.

Another class of approaches includes delayed acceptance MCMC algorithms, which utilize cheaper model approximations to accelerate the parameter search (also via proposal distribution modification)~\cite{christen2005markov}. 
Several works have been completed that develop and explore the use of cheaper models with coarser spatial discretizations in a two-stage or multilevel framework. Early works that employ coarser spatial discretizations include \cite{higdon2002bayesian} and \cite{efendiev2006preconditioning}; the former utilizes a Metropolis coupled MCMC to swap proposals between coarse and fine chains, and the latter performs a delayed acceptance where sample proposals are only completed with the fine grid solver if their associated coarse grid solutions have been accepted. More recent works have investigated multilevel MCMC approaches. In particular, in \cite{Dodwell15}, the authors developed an approach to both accelerate the mixing of the MCMC chain by using multiple levels, each with coarser spatial discretizations, and accelerate the sampling by performing variance reduction via multilevel Monte Carlo following the ideas of ~\cite{Heinrich01, Giles08, barth2011multi, Cliffe11, Teckentrup13}.
Analysis of a multilevel MCMC was completed in~\cite{hoang2019analysis}. While promising speed up results have been shown, numerical testing has been limited to 2D spatial domains \hrev{and structured meshes}.

\hrev{While scalable solvers are available for several classes of PDE forward models, the sampling of large-scale Gaussian random field on unstructured meshes in an algorithmically scalable manner is still a challenging task. The use of a Karhunen-Lo\`{e}ve (KL) expansion to form Gaussian random field realizations requires calculating the eigenvalues and eigenfunctions of the covariance function~\cite{ghanem2003stochastic}. A straightforward, though perhaps na\"{i}ve
implementation will have a cost that grows cubically with the degrees of freedom associated with the spatial discretization of the random field, i.e., the mesh size. While there are tools to improve this scaling, e.g., hierarchical matrix formations \cite{bebendorf2008hierarchical} or Nystr\"{o}m methods~\cite{williams2001using}, storage and the ability to calculate the KL expansions for unstructured meshes are roadblocks to large-scale and extreme-scale applications. Other approaches to sampling, such as circulant embedding~\cite{graham2018analysis}, are not directly applicable to problems with unstructured meshes. }

An alternative scalable approach to generate random field realizations is via the stochastic reaction-diffusion PDE formulated in \cite{Whittle63} and solved with finite elements in \cite{Lindgren11}. Using this approach, each independent realization requires solving the stochastic PDE with an independent realization of spatial white noise function as forcing term. Applying this approach in a multilevel setting, such as multilevel Monte Carlo or multilevel MCMC requires forming coupled realizations of Gaussian random fields on multiple levels of discretization. A few works have completed this, including \cite{drzisga2017scheduling} where fine and coarse level realizations are coupled together to perform massively parallel multilevel Monte Carlo. In \cite{Osborn17, Osborn17b} the authors solve a mixed PDE on the space of piecewise constants, and generate matching fine and coarse realizations by restricting the fine grid spatial white noise to the coarse level, using operators and solvers from element agglomerated algebraic multigrid (AMGe). In \cite{croci2018efficient} the authors couple the coarse and fine level realizations using the primal formulation of the PDE.

While these approaches have been incorporated successfully into the multilevel Monte Carlo framework, they are not useful in the multilevel MCMC framework. \rsout{This is because in the aforementioned stochastic PDE approaches, the fine and coarse level realizations are coupled across the two levels. For example, in \cite{Osborn17, Osborn17b}, the fine level realizations are formed by sampling the spatial white noise source function on the fine level (to form the fine level Gaussian random field realization), and then the corresponding coarse realizations are formed by restricting the fine spatial white noise source function to the coarse level (to form the fine level Gaussian random field realization).}\hrev{This is because, in the multilevel MCMC approach, }we must first sample from the coarse level, and then form a fine level random field in a hierarchical manner from the coarse realization. As this sampling approach has not yet been developed (to the best of the authors' knowledge), this paper seeks to fill this void.

\subsection{Contributions of this Work}
In this work, we develop \hrev{an algorithmically} scalable, hierarchical Gaussian random field sampling method that can be used to construct proposals distributions in the multilevel MCMC framework. Specifically, we plug our sampling method into the multilevel MCMC framework of \cite{Dodwell15}, though it is also applicable to other two-level MCMC or delayed acceptance MCMC approaches discussed earlier. To do this, we utilize the finite element solvers from \cite{Osborn17b} to \hrev{map an independent realization of spatial white noise to a Gaussian random field realization. Within this mapping, we incorporate a new component: a hierarchical decomposition of the white noise (via $L^2$ projection operators) across discretization levels. This new feature allows us to perform MCMC stepping on coarse level white noise, extend it to a finer level, and then perform an MCMC step on independent white noise in the complementary space.}

The remainder of this paper is organized as follows. In Section \ref{sec:grf}, mathematical notation relevant to Gaussian random fields is presented\rsout{ to provide a framework for this sampling approach. In addition this section presents the stochastic PDE approach to calculating discrete Gaussian random field realizations}.
Our new hierarchical approach is presented in Section \ref{sec:hspde}; this includes the theoretical aspects of performing a hierarchical direct decomposition of white noise -- in a two-level and multilevel framework -- resulting in a hierarchical approach to form Gaussian random field realizations.\rsout{These aspects are presented in Lemma \ref{thm:2L} and Theorem \ref{thm:ML}.} The numerical implementation is discussed in Section \ref{sec:hspde_implementation}, in the form of algorithms, as well as visualizations of the random field hierarchical decomposition. \rsout{Numerical results are provided in Sections \ref{sec:num1} and \ref{sec:num2}.} Section \ref{sec:num1} explores the cost and scaling of our multilevel hierarchical sampling technique applied to the Egg Model \cite{jansen2014egg} using three levels\hrev{; in particular, we show that the algorithm is scalable}. Section \ref{sec:num2} incorporates this new hierarchical sampling technique into a \hrev{four}-level MCMC following the approach of~\cite{Dodwell15}, and shows that we obtain similar improvements in the multilevel acceptance rate, variance decay, and total computational cost when compared to the single-level approach.

\subsection{Mathematical Notation}
As a reference to the reader, we define the majority of this paper's notation in Table \ref{tab:notation1}. The first section of the table introduces general variables that provide a basis for the majority of the mathematical notation. The second section of the table refers to discrete variables that are used in various finite element representations, and that are frequently referred to throughout this work.

\begin{table} \footnotesize
\caption{Mathematical Notation.}\label{tab:notation1}
\begin{center}\begin{tabular}{ll}
\toprule 
Variable   & Description\\                              \midrule
 $\x \in D \subset \mathbb{R}^d$ & Point in spatial domain, $d=2 $ or $3$  \\%
 $\omega \in \Omega$ & Outcome of Sample Space \\
 $ u \in \Theta := L^2(D)$ &  $L^2$ function defined over $D$\\
 $\zeta:=\zeta(\x , \omega)$ & White noise function in $D$\\
 $ \bq \in \bm R := H(\text{div}; D) $ & $H(\text{div})$ function defined over $D$\\
 \midrule
 Discrete Variable & \\
 \midrule
 $h,H $ & Subscripts to denote fine and coarse level objects \\
 $\ell$ & Subscript denoting running level $\ell$ index, with $\ell =0$ as finest \\
 $k$ & Subscript denoting target level of an algorithm \\
 $\TT _\ell$ &  Level $\ell$ finite element triangulation \\
 $u_\ell\in \Theta_\ell$ & Piecewise constant function defined on $\TT_\ell$ \\
 $\Q _\ell$ & Orthogonal projection from $L^2$ to $\Theta _\ell$ \\ 
 $P_\ell$ & Interpolation operator mapping between $\Theta _{\ell +1}$ and $\Theta_\ell$  \\
 $\Pi _\ell$ & Restriction operator mapping between $\Theta _{\ell }$ and $\Theta_{\ell+1}$ \\
  $ \zeta _\ell$ & White noise representation in $\Theta_\ell$ \\
  $\bzeta _\ell$ & Coefficient vector of white noise finite element representation in $\Theta_\ell$\\
  $\bxi _\ell$ & Vector of random elements \\
  $\bb _\ell$ & Vector representation of white noise in $\Theta _\ell$ \\
 $\bq _\ell\in \bm R _\ell$ & Function of the lowest order Raviat-Thomas space on $\TT_\ell$ \\
    $M_\ell, B_\ell, W _\ell$ & Mass matrices for various level $\ell$ spaces \\
   \bottomrule                                                                         
    \end{tabular}\end{center}  
\end{table}  

\section{Gaussian Random Fields}\label{sec:grf}
In this work we consider a particular class of random fields, that is, spatially correlated Gaussian random fields, which in this context, will be used to describe an uncertain physical process. Define the probability space $(\Omega, \mathcal{F}, \mathbb{P} )$, with sample space $\Omega$, $\sigma$-algebra $\mathcal{F}$, and probability $\mathbb{P} $. Given the spatial domain of interest $D\subset \R^d$, with $d = 2, 3$, we seek to form random field realizations of $\{u( \x, \omega)\in L^2(D): \x \in D,\ \omega \in \Omega \}$, that follow a Gaussian prior density $u \sim \mathcal{N}(0, \mathcal{C})$ with zero mean and covariance operator $\mathcal{C}$. To ensure the mesh independent statistics of the random field $u$, $\mathcal{C}$ is a trace-class operator \cite{Stuart2010}. Specifically, we define the covariance operator as the squared inverse elliptic operator (see e.g. \cite{Flath11, bui2013computational, petra2014computational}). That is, 
\begin{equation}\label{eq:cov_op}
\mathcal{C} = \mathcal{A}^{-2} \text{ with } \mathcal{A} u := -\div \left( \frac{1}{g} \nabla u\right) + \frac{\kappa^2}{g} u,
\end{equation}
where  $\kappa$ denotes the inverse of the correlation length and $g$ controls the marginal variance of the field.
Using the above notation, we then define the probability density function as
\begin{equation}\label{eq:prior}
d\mu(u) \propto \exp\left( -\frac{1}{2} \langle \mathcal{A} u,  \mathcal{A} u\rangle \right),
\end{equation}
where $ \langle \mathcal{A} u,  \mathcal{A} u\rangle = \int_D (\mathcal{A} u)^2 d\mathbf{x}$.

As described in \cite{Whittle63,Lindgren11}, for unbounded domains $D:= \mathbb{R}^d$, the covariance operator in \eqref{eq:cov_op} leads to a Gaussian random field of the Mat\'ern family with smoothness parameter $\nu$ and marginal variance $\sigma^2$ respectively given by
$$ \nu = 2 - \frac{d}{2} \text{ and } \sigma^2 = \frac{ g^2 \Gamma(\nu)}{\Gamma(\nu + d/2) (4\pi)^{d/2}\kappa^{2\nu}}.$$
In particular, in three-spatial dimensions, this gives the well-known exponential covariance operator

\begin{equation*}  
(\mathcal{C}u)(\mathbf{x}) := \int_D \text{cov}(\mathbf{x},\mathbf{y}) u(\mathbf{y}) d\mathbf{y}, \text{ with } \text{cov}(\mathbf{x},\mathbf{y}) :=\frac{g^2}{8\pi\kappa} \exp \left(-\kappa \norm{\x - \y}_2 \right).
\end{equation*}

For a finite domain $D \subset \mathbb{R}^d$, suitable boundary conditions need be stipulated to reduce boundary artifacts, see e.g. \cite{roininen2014whittle,khristenko2019analysis,daon2018mitigating}. In this work, we choose to extend the domain $D$ to a larger domain $\overline{D}  \subset \mathbb{R}^d$ and equip $\mathcal{A}$ with homogeneous Neumann boundary conditions on $\partial{\overline{D}}$, as described in \cite{Osborn17b}.

For sample-based UQ approaches---such as standard Monte Carlo---we desire to generate samples of this random field $u( \x, \omega)$ to serve as input field data to a model of interest. In our application (which we further detail in Section \ref{sec:num2}), we wish to generate permeability field realizations, $k = \exp(u( \x, \omega))$, each of which serves as an input to Darcy's equations.

\subsection{A Stochastic PDE Approach for Finite Element Random Fields}\label{sec:spde}
As presented in \cite{Flath11, bui2013computational, petra2014computational}, a realization of a Gaussian random field, with covariance operator $\mathcal{C}$ given by \eqref{eq:cov_op}, can be generated by solving the stochastic reaction-diffusion PDE
$$ \mathcal{A} u =\zeta, $$
where $\zeta := \zeta(\bf x, \omega)$ is spatial Gaussian white noise. The spatial Gaussian white noise $\zeta$ is an $L^2(D)$-bounded generalized function \cite[Appendix B]{Lindgren11}, such that
\begin{equation}\label{eq:whiteNoiseSourceFunctionDef}
 \left <\zeta, v\right> \sim \mathcal{N}(0, \|v\|^2_{L_2(D)}) \quad \forall v \in L^2(D).
 \end{equation}
In the following, we consider a particular PDE-based approach that uses a mixed formulation to generate field realizations. That is, we follow the approach of \cite{Osborn17, Osborn17b}, which allows us to work in the space of piecewise constants. For large-scale applications this is beneficial as it provides a natural way to define spatial white noise, and the associated mass matrix is easily diagonalizable. 

\subsubsection{{A Mixed Formulation}}
\rsout{The following mixed stochastic PDE used in \cite{Osborn17, Osborn17b} provides a way to generate Gaussian random field realizations:}
For a fixed $\omega \in\Omega$, a Gaussian random field realization $u := u(\x, \omega)$ is calculated by solving the stochastic PDE:
    \begin{equation}
       \begin{array}{ll}
            (\brho, \s)+ (\div \s, u)= 0 & \forall \s \in H(\div) \\
            (\div \brho, v)- \kappa^2\, (u, v) = -g\, \left <\zeta,v\right >& \forall v \in L^2,
        \end{array}\label{eq:spde_weak1}
    \end{equation}
where $(\cdot, \cdot)$ denotes the $L^2(D)$ inner product~\cite{Osborn17, Osborn17b}.
Above, the spatial Gaussian white noise $\zeta := \zeta(\bf x, \omega)$ is a zero-mean random Gaussian field on $D$ such that $\left <\zeta,v\right > \sim \mathcal{N}(0, \| v \|^2_{L^2(D)})$, for any function $v \in L^2(D)$ (see \eqref{eq:whiteNoiseSourceFunctionDef}). Note, properties of finite element white noise will be discussed in the following section.

Define the spaces $\Theta = L^2(D)$ with inner product $(u , v) = \int_D u v d\x$ for all $u, v \in \Theta$ and $\bm R = H(\text{div}; D):= \{ \bq \in \bm [L^2(D)]^d | \text{ div }\bq \in L^2(D),\ \bq \cdot \bn =0 \text{ on }\partial D \}$ with inner product $(\bq , \s)= \int_D \bq \cdot \s d\x$ for all $\bq, \s \in \bm R$. Let $\bR_h, \Theta_h$ be the pair of the lowest order Raviart-Thomas and piecewise constant finite element spaces associated with the given triangulation $\T_h$.

For a fixed $\omega \in \Omega$, discrete solutions $\brho_h \in \bR_h$ and $u_h \in \Theta_h$ are calculated from the mixed system,
    \begin{equation}
       \begin{array}{ll}
            (\brho_h, \s_h)+ (\div \s_h, u_h)= 0 & \forall \s_h \in \mathbf{R}_h \\
            (\div \brho_h, v_h)- \kappa^2\, (u_h, v_h) = -g\, \left < \zeta,v_h \right > & \forall v_h \in \Theta_h.
        \end{array}\label{eq:spde_weak}
    \end{equation}

\subsubsection{Finite Element Representation of White Noise}
\urev{Since moments of $\zeta$ are well-defined for functions in $\Theta_h$, we can define the mapping $\Q_h: \zeta \mapsto \Q_h \zeta \in \Theta_h$ using the identity
\begin{equation}\label{eq:wn_eqiv}
(\Q_h \zeta,v_h )  =\left <\zeta,\; v_h\right > \ \forall v_h\in \Theta_h.
\end{equation}
That is, a realization of white noise on a given finite element mesh $\mathcal{T}_h$ can be represented in $\Theta_h$ using the mapping $\Q_h$ as follows
\begin{equation}\label{eq:wn_exp}
\zeta_h:= \Q_h \zeta = \sum_{\tau \in \T_h} \zeta_\tau \chi_\tau ,
\end{equation}
where $\{\chi_\tau\}$ is an $L^2$-orthogonal basis of piecewise constants spanning $\Theta_h$.}

\urev{Using the expansion in \eqref{eq:wn_exp} and the equivalence in \eqref{eq:wn_eqiv}, it follows that the righthand side of \eqref{eq:spde_weak} will have the coefficient vector $\bb_h \equiv\left (
( \zetacoeff_h,\;\chi_\tau) \right )_{\tau \in \T_h}$. As a consequence of using piecewise constant basis functions, each inner product simplifies as
\begin{equation*}
\bb_h = ( \zetacoeff_h,\;\chi_\tau) =  \zetacoeff_\tau \|\chi_\tau\|^2.
\end{equation*}
In other words, $\bb_h = \bW_h \bzeta_h$, where $\bW_h$ the diagonal mass matrix for the space $\Theta_h$ and $\bb_h$ is the vector collecting the coefficients $\zeta_\tau$ in the expansion \eqref{eq:wn_eqiv}.}

To generate realizations of white noise in $\Theta_h$, we consider the following properties (see \cite[Section 1.4.3]{adler2009random} and \cite[Section 2.4.5]{Adler07applicationsof} for details). \hrev{We note that, while we present these properties with respect to an $L^2$-orthogonal basis of piecewise constants, they can be generalized to the situation of a non-orthogonal basis.}

\begin{property}[White noise in $\Theta _h$]\label{def:white_noise_h}
Let $\zeta$ be white noise in $D$. Then, for the projection of $\zeta$ onto the basis $\{\chi_\tau\}_{\tau\in\TT_h}$ of $\Theta_h$, denoted $ \zetacoeff_h $ as in (\ref{eq:wn_exp}), it follows that,
\begin{equation*}
 \E[( \zetacoeff_h,\;\chi_\tau) ] =\E[\left < \zeta,\;\chi_\tau \right >] = 0,
\end{equation*}
and
\begin{equation*}
 \E[( \zetacoeff_h,\;\chi_{\tau_i}) ( \zetacoeff_h,\;\chi_{\tau_j})] = \E[\left < \zeta,\;\chi_{\tau_i} \right > \left <\zeta,\;\chi_{\tau_j} \right >] =
(\chi_{\tau_i},\;\chi_{\tau_j}),
\end{equation*} 
which implies 
\begin{equation*}
\E[\left (( \zetacoeff_h,\;\chi_{\tau_i})\right )_{\tau_i\in\TT_h}  (  \left ( ( \zetacoeff_h,\;\chi_{\tau_i})    \right ) _{\tau_j\in\TT_h}  )^T] 
= \left ((\chi_{\tau_i} ,\;\chi_{\tau_i}) \right  )_{\tau_i,\tau_j \in\TT_h} = \bW_h,
\end{equation*}
where $\bW_h$ is the (diagonal) mass matrix for the space $\Theta_h$. 
\end{property}

These properties follow from the theoretical aspects of white noise. Specifically, the covariance between two volumes $A$ and $B$ (within $D$) is equivalent to the mass of the intersection of the two volumes (further theoretical aspects of Gaussian white noise may be found in \cite{adler2009random,Adler07applicationsof}), and for finite element white noise this implies that the covariance is equivalent to the mass matrix.
\hrev{Using the above properties, we can show that for $\bW_h \bzeta_h$ to be a realization of Gaussian white noise, we require $\bzeta_h= (\zeta_\tau)_{\tau \in \T_h}\sim \N(0,\;\bW^{-1}_h)$.}
\rsout{$\bzeta_h = ( \zetacoeff_i)^n_{i=1} \sim \N(0,\;\bW^{-1}_h)$}
\begin{lemma}\label{lem:wn}
Given the basis \hrev{$\{\chi_\tau\}_{\tau\in\TT_h}$ of $\Theta_h$, associated mass matrix $\bW_h =\left ((\chi_{\tau_i} ,\;\chi_{\tau_i}) \right  )_{\tau_i,\tau_j \in\TT_h}$, and $\bzeta_h= (\zeta_\tau)_{\tau \in \T_h}$} sampled from $\N(0,\;\bW^{-1}_h)$, it follows that $\bW_h \bzeta_h$ is a realization of white noise in $\Theta_h$.
\end{lemma}

\begin{proof}
Following Property \ref{def:white_noise_h}, it is sufficient to show that $\E[\bW_h \bzeta_h]=0$ and $\E[\bW_h \bzeta_h(\bW_h \bzeta_h)^T]=\bW_h$. As $\bzeta_h \sim \N(0,\;\bW^{-1}_h)$, it is clear that $\E[\bW_h \bzeta_h]=0$. As for the covariance, we have 
\begin{equation*}
\begin{array}{lcl}
\E[\bW_h \bzeta_h(\bW_h \bzeta_h)^T] &= & \bW_h  \E[\bzeta_h\bzeta_h^T]\bW_h \\
&= & \bW_h .
\end{array}
\end{equation*}

\end{proof}

\subsubsection{Finite Element Representation of Gaussian Random Fields}
\label{sec:fe_grf}
In the actual computation of $u_h$, we use the equivalent vector representation for the righthand side of \eqref{eq:spde_weak}, defined as
\begin{equation}\label{eq:bh}
-g \bb_h = -g \bW^{1/2}_h { \bxi}_h,
\end{equation}
with ${ \bxi}_h\sim~\N(0, \;I)$. We note this equivalence is made clearer in Section \ref{sec:hspde_implementation}. As we are in the space of piecewise constants in $L^2$, the square root of the \hrev{(diagonal)} mass matrix $\bW_h$ is easily calculated. Let $M_h$ be the mass matrix associated with inner product $(\brho_h, \s_h)$ and $B_h$ the mass matrix associated with the bilinear form $(\div \s_h, u_h)$.
Then the matrix representation of \eqref{eq:spde_weak} is given as
\begin{equation}\label{eq:mixed_system}
\begin{bmatrix} M_h & B_h^T \\ B_h & -\kappa^2 W_h \end{bmatrix} \begin{bmatrix} \brho_h \\ \bu_h \end{bmatrix} = 
\begin{bmatrix} \boldsymbol{0} \\ -g \bb_h \end{bmatrix},
\end{equation}
with $\bb_h$ defined by \eqref{eq:bh}.

For ease of notation, we introduce the scaled negative Schur Complement of \eqref{eq:mixed_system} defined by
\begin{equation}
A_h :=  \frac{\kappa^2}{g} W_h + \frac{1}{g} B_h M_h^{-1} B_h^T.
\end{equation}
As demonstrated in \cite{Osborn17}, solutions $u_h$ of the mixed system in \eqref{eq:mixed_system} are discrete realizations of a Gaussian random field with density $\mu_h \sim \mathcal{N}( \boldsymbol{0}, C_h)$, where $C_h = A_h^{-1} W_h  A_h^{-1}$. It then follows that the corresponding probability density function is
\begin{equation}
\mu_h(u_h) \propto \exp{ (- \bu_h^T A_h W_h^{-1} A_h \bu_h)} = \exp{ (- \bb_h^T W_h^{-1} \bb_h)}.
\end{equation}

\section{Multilevel Hierarchical Decomposition of Finite Element White Noise}\label{sec:hspde}

In what follows we study the computational aspects of sampling the righthand side in \eqref{eq:spde_weak} from a coarse finite element space $\Theta_H \subset \Theta_h$, and its (direct) hierarchical complement space 
$(I-\Q_H)\Theta_h$, where $\Q_H:\; L^2 \mapsto \Theta_H$ is the corresponding $L^2$-projection. For any $ \zetacoeff_h \in \Theta_h$, we use the two-level hierarchical decomposition
\begin{equation*}
 \zetacoeff_h = \Q_H  \zetacoeff_h + (I-\Q_H) \zetacoeff_h
\end{equation*}
to decompose $ \zetacoeff_h$ into the spaces $\Theta _H$ and $\Theta _h\backslash \Theta_H$. Since we work with spaces of discontinuous (piecewise constant) functions $\Theta_h$ and $\Theta_H$ with associated mass matrices $\bW_h$ and $\bW_H$, respectively, the projection $\Q_H$ is easily implemented (by inverting a diagonal (coarse) mass matrix). 

Define $P$ to be the interpolation matrix that relates the coarse coefficient vector $\bzeta_H$ of $ \zetacoeff_H$ (expanded in terms of the basis of $\Theta_H$) and the fine coefficient vector $\bzeta_h$ of $ \zetacoeff_H \in \Theta_H \subset \Theta_h$ expanded in terms of the basis of $\Theta_h$. That is, 
\begin{equation*}
\bzeta_h = P \bzeta_H.
\end{equation*}
Let $\bPi = W_H^{-1} P^T W_h$ denote the restriction operator, then $P\bPi$ is the matrix representation of $\Q_{H}$ and $\bPi P = I$. That is, we have $P \bzeta_H =  P\bPi \bzeta_h$, or 
$\bzeta_H = \bPi \bzeta_h$. 

In what follows, we first seek to show that $ \Q_H  \zetacoeff_h$ gives rise to a coarse random coefficient vector $\bzeta_H \sim \N(0,\;\bW^{-1}_H)$.
\begin{lemma}\label{lem:2L_coarse}
Let $ \zetacoeff_h \in \Theta_h$, with coefficient vector $\bzeta_h \sim\N(0,\;\bW^{-1}_h)$. Then $ \zetacoeff_H \equiv \Q_H  \zetacoeff_h$ has coefficient vector $\bzeta_H \sim \N(0,\;\bW^{-1}_H)$.
\end{lemma}
\begin{proof}
Given the associated coarse coefficient $\bzeta_H = \bPi \bzeta_h$ with $\bzeta_h \sim \N(0,\;\bW^{-1}_h)$,  it is clear the mean is zero. For the covariance matrix, we have
\begin{equation*}
\begin{array}{rl}
\E[\bzeta_H\bzeta^T_H] &= \E[\bPi\bzeta_h (\bPi \bzeta_h)^T]\\
&=\bW^{-1}_H P^T \bW_h \E[\bzeta_h \bzeta^T_h] \bW_h P \bW^{-1}_H\\
& = \bW^{-1}_H.
\end{array}
\end{equation*}
Above, we use $ \E[\bzeta_h \bzeta^T_h] = \bW^{-1}_h$ and the Galerkin relation between the coarse and the fine level mass matrices,
$\bW_H = P^T \bW_h P$. 
\end{proof}
\rsout{Therefore, if $ \zetacoeff_h$ has coefficient vector $\bzeta_h \sim \N(0,\;\bW^{-1}_h)$ then the same holds for the coarse projection $ \zetacoeff_H = \Q_H  \zetacoeff_h$, that is, its coefficient vector $\bzeta_H = \bPi \bzeta_h \sim \N(0,\;\bW^{-1}_H)$. If $ \zetacoeff_h$ is a fine finite element representation of white noise, then 
its projection $ \zetacoeff_H = \Q_H  \zetacoeff_h$ is the corresponding coarse finite element representation of white noise. }

Next we present our main lemma, which allows us to utilize this two-level, hierarchical decomposition to form a realization of white noise on $\Theta_h$. 
\begin{lemma}\label{thm:2L}
Let $ \zetacoeff_H \in \Theta_H$ be a coarse representation of white noise with a coarse coefficient vector $\bzeta_H \sim \N(0,\; \bW^{-1}_H)$, and let $ \zetacoeff_h \in\Theta_h$ be a fine representation of white noise with fine coefficient vector $\bzeta_h \sim \N(0,\; \bW^{-1}_h)$, such that $\bzeta_H$ and $\bzeta_h$ are independent. Then the fine level function
\begin{equation}\label{eq:hdecomp}
 \zetacoeff^{'}_h =  \zetacoeff_H + (I-\Q_H)  \zetacoeff_h
\end{equation}
is a representation of the white noise in $\Theta_h$.
\end{lemma}
\begin{proof}
First, consider the coefficient vector of $ \zetacoeff^{'}_h$, given as
\begin{equation}\label{eq:2levelWN}
\bzeta^{'}_h = P \bzeta_H +(I - P\bPi) \bzeta_h.
\end{equation}
To prove $ \zetacoeff^{'}_h$ is a representation of Gaussian white noise, we must show Definition \ref{def:white_noise_h} holds, that is $\E[\bzeta^{'}_h]=0$, and $\E[\bzeta^{'}_h (\bzeta^{'}_h)^T]= \bW_h^{-1}$. We assume  that $\bzeta_H$ and $\bzeta_h$ are independent, which  implies that 
$\E[\bzeta_H \bzeta^T_h] = \E[\bzeta_H] \E[\bzeta^T_h] = 0$. 
Hence for the covariance matrix, we have
\begin{equation*}
\begin{array}{rl}
\E[\bzeta^{'}_h \left (\bzeta^{'}_h \right )^T]
& \displaystyle = \E [ (P \bzeta_H +(I - P\bW^{-1}_H P^T \bW_h) \bzeta_h) (P \bzeta_H +(I - P\bW^{-1}_H P^T \bW_h) \bzeta_h)^T]\\
&\displaystyle = P\E[\bzeta_H \bzeta^T_H]P^T + (I - P\bW^{-1}_H P^T \bW_h) \E[\bzeta_h\bzeta^T_h] (I - \bW_h P \bW^{-1}_H P^T)\\
& \displaystyle = P \bW^{-1}_H P^T +(I - P\bW^{-1}_H P^T \bW_h) \bW^{-1}_h(I - \bW_h P \bW^{-1}_H P^T)\\
&\displaystyle = \bW^{-1}_h.
\end{array}
\end{equation*}
That is, $\bzeta^{'}_h \sim \N(0,\;\bW^{-1}_h)$; hence $ \zetacoeff^{'}_h$ is a fine finite element representation of white noise.
It is clear also that $ \zetacoeff_H = \Q_H  \zetacoeff{'}_h$ and $(I-\Q_H)  \zetacoeff_h = (I-\Q_H)  \zetacoeff{'}_h$. 
\end{proof}

In conclusion, the finite element hierarchical (direct) decomposition based on $\Q_H$ provides a hierarchical decomposition of the fine finite element white noise 
into a coarse finite element representation of white noise plus a computational hierarchical (direct) complement which also involves 
fine finite element representation of white noise. 

\subsection{The multilevel hierarchical decomposition}
\hrev{To extend the the above two-level hierarchical decomposition of Gaussian white noise to a multilevel hierarchical decomposition, we introduce the following notation.}
\rsout{Next we wish to extend the above two-level hierarchical decomposition of Gaussian white noise to a multilevel hierarchical decomposition, so that we may sample across multiple levels of discretization.}
Let $\TT_0\equiv \TT_h$ denote the finest level triangulation of $D$, with a hierarchy of $L$ coarser levels given as $\{\TT_\ell \}_{\ell=1}^L$, such that $\TT_L$ represents the coarsest triangulation. We consider the finite element space $\Theta _\ell $ to be the space of piecewise constant functions associated with the triangulation $\TT_\ell$, for $\ell = 0, \hdots, L$, and with mass matrix $\bW_\ell$; and $\bR _\ell $ the lowest order Raviart-Thomas space associated with the triangulation $\TT_\ell$. Additionally, define the sequence of $L^2$-projections $\Q_\ell:\; L^2 \mapsto \Theta_\ell$ with $\ell = 0,\hdots L$. 

In what follows, we construct the multilevel hierarchical decomposition of white noise for a given level $k<L$.
\begin{thm}\label{thm:ML}
Consider the representations of white noise, given as $ \zetacoeff_\ell \in \Theta_\ell$ with associated coefficient vectors $\bzeta_\ell \sim \N(0,\; \bW^{-1}_\ell)$, for $\ell = k,\hdots, L$, such that each $\bzeta_\ell $ is independent. Then the level $k$ function
\begin{equation}\label{eq:mlhdecomp}
 \zetacoeff^{'}_k =  \zetacoeff_L + \sum\limits_{\ell = k}^{L-1} (I-\Q_{\ell+1})  \zetacoeff_{\ell},
\end{equation}
with $k<L$, is a representation of the white noise in $\Theta_k$.
\end{thm}

\begin{proof}
\hrev{From Lemma \ref{lem:2L_coarse} the result is clear for $k=L-1$, i.e., the two-level case. For additional levels, the result follows by applying Lemma \ref{lem:2L_coarse} in a recursive manner.}
\rsout{First, consider the two-level decomposition. Let $ \zetacoeff_L \in \Theta_L$ be a realization of white noise from the coarsest level (level $L$). Given an independent realization of white noise from level $L-1$, denoted $\zetacoeff_{L-1} \in \Theta_{L-1}$, it follows from from Lemma \ref{thm:2L} that 
\[
 \zetacoeff^{'}_{L-1} =  \zetacoeff_{L} + (I-\Q_{L})  \zetacoeff_{L-1},
\]
is a realization of white noise in $\Theta_{L-1}$. To generate a realization on level $L-2$, we simply sample independent white noise $\zetacoeff^{'}_{L-2} \in \Theta_{L-2}$, and form the hierarchical realization that builds from our previous white noise realization $\zetacoeff^{'}_{L-1}$:     
\[
 \zetacoeff^{''}_{L-2} =  \zetacoeff^{'}_{L-1} + (I-\Q_{L-1})  \zetacoeff^{'}_{L-2},
\]
which, by Lemma \ref{lem:2L_coarse}, is a realization of white noise in $\Theta_{L-2}$. Continuing in this fashion to level $k$, let $ \zetacoeff^{*}_{k} \in \Theta_k$ be a sample of white noise that has been hierarchically formed. Let $ \zetacoeff^{*}_{k-1} \in \Theta_{k-1}$ be an independently sampled realization of level $k-1$ white noise. Then it follows from Lemma \ref{lem:2L_coarse} that
\[
 \zetacoeff^{**}_{k-1} =  \zetacoeff^{*}_{k-2} + (I-\Q_{k-2})  \zetacoeff^*_{k-1},
\]
is a realization of white noise in $\Theta_{k-1}$. By induction, it follows that, for an arbitrary level $\ell$, this hierarchical multilevel construction formulates a realization of white noise in $\Theta_\ell$.}
\end{proof}

The associated coefficient representation is defined similarly to (\ref{eq:2levelWN}); however, we replace the subscript $h$ with $k$ to denote the level, and let $P_k$ be the interpolation matrix that maps the level $k+1$ coefficient vector $\bzeta_{k+1}$ of $ \zetacoeff_{k+1}$ to the
level $k$ coefficient vector $\bzeta_{k}$ of $ \zetacoeff_{k}$, such that $\bzeta_k = P_k \bzeta_{k+1}$. This hierarchical coefficient representation (in two-level form) is given as
\begin{equation}\label{eq:2levelWN_k}
\bzeta^{'}_k = P_k \bzeta_{k+1} +(I - P_k\bW^{-1}_{k+1} P_k^T \bW_{k}) \bzeta_k.
\end{equation}
Just as in the proof, we can hierarchically build a level $k$ coefficient vector by starting on the coarsest level and adding on coefficients projected onto the complimentary spaces, as will be further detailed in the next section (see, e.g., Algorithm \ref{alg:hspde}).

\rsout{The multilevel, hierarchical decomposition in (\ref{eq:mlhdecomp}) provides the theoretical framework that allows us to decompose the white noise across multiple levels, enabling us to perform hierarchical sampling of Gaussian random field.}

\section{Implementation of the Hierarchical Sampler}\label{sec:hspde_implementation}
\urev{Recall from Lemma \ref{lem:wn} that we may sample (single-level) finite element white noise on level $k$ via $b_k = \bW_k \bzeta_k$ with $\bzeta_k = ( \zetacoeff_i)^n_{i=1} \sim \N(0,\;\bW^{-1}_k)$. 
While we can use the decomposition in (\ref{eq:2levelWN_k}) for our hierarchical implementation, we instead \hrev{alter} this representation to accomplish two additional goals: first, that on each level we sample from a $\N(0,I)$ distribution, and second, that we utilize the interpolation and restriction operators $P_k$ and $\bPi_k=W_{k+1}^{-1} P_k^TW_k$.} %

\urev{By multiplying \eqref{eq:2levelWN_k} by $\bW_k$, and after simple algebraic manipulation, we obtain the following hierarchical representation of the right hand side  of (\ref{eq:spde_weak}):
\begin{equation*}
\bb_k^{'}  = \bPi^T_k \bb_{k+1} + (I - \bPi^T_kP_k^T) \bb_k,
\end{equation*}
where $\bb_k^{'} =  \bW_k { \bzeta}_k^{'}$.}

For algorithmic efficiency (and to meet our additional two goals), we simplify the above using the fact that $\bb_\ell = \bW^{1/2}_\ell { \bxi}_\ell$ with ${ \bxi}_\ell \sim~\N(0, \;I)$ for $\ell = k, k+1$ and write
\begin{equation*}
b_k^{'} = \bW^{1/2}_k{ \bxi}_k^{'}   =\bPi^T_k (\bW^{1/2}_{k+1} { \bxi}_{k+1}) + (I -\bPi^T_kP_k^T)(\bW^{1/2}_k { \bxi}_k).
\end{equation*} 
In practice, \hrev{we use this finite element white noise formulation of $\bb_k$, where we may construct a realization using multiple coarser levels, beyond that of level $k+1$.}
This process is described in Algorithm \ref{alg:hspde}. For a given level $k$, we first calculate finite element white noise on the coarsest level $L$, denoted $\bb_L$. Then we iterate through each finer level, where we calculate $\bb_\ell$ by first interpolating the coarser $\bb_{\ell+1}$ (which was previously calculated), and then adding a spatial white noise realization that is complementary to the coarser $\ell+1$ space -- this is accomplished by multiplying level $\ell$ spatial white noise, $\bW_\ell^{1/2}\bxi_\ell$, with $(I - \bPi^T_\ell P^T_\ell)$, which projects the level $\ell$ spatial white noise orthogonal to the coarser space(s). After each iterate, we refine a level (decrease $\ell$ by 1), and repeat this process. This done until we reach level $k$, and the resulting $\bb_k$ provides us with our hierarchically generated realization of spatial white noise, which can then be used in the righthand side of the discrete problem \eqref{eq:spde_weak}.

\RestyleAlgo{boxruled}
\SetKwInOut{Input}{input}
\SetKwInOut{Output}{output}
\SetKw{Initialize}{initialize}
\SetKw{Return}{return}
\SetKw{Break}{break}
\begin{algorithm}
\KwIn{Current level $k$ (with $0 \leq k \leq L$), $L$, independent $\{\bxi_L, \bxi_{L-1},\hdots ,\bxi_k \}$ with each $\N(\bm 0, I_k$) }
$\ell = L$\\
$\bb_L = \bW_L^{1/2}\bxi_L$\\
$\ell = \ell-1$\\
\While{$\ell\geq k$}{
$\bb_\ell = \bPi^T_\ell \bb_{\ell+1} + (I - \bPi^T_\ell P^T_\ell) \bW_\ell^{1/2}\bxi_\ell$ \\
$\ell= \ell-1$
}
\KwOut{$\{\bb_{L}, \bb_{L-1}, \hdots , \bb_k \}$}
\protect\caption{Form finite element white noise via new hierarchical approach. %
\label{alg:hspde}}
\end{algorithm}

In this work we employ the same linear system of \cite[Section 2.2]{Osborn17}, but instead of spatial white noise generated strictly on the fine level, we use our hierarchical approach. For a given level $k$, we seek to calculate solutions $(\brho_k,u_k ) \in \bR_k \times \Theta_k $ via the linear system%
\begin{equation} 
    \begin{array}{lcr}
        \begin{bmatrix} 
            M_k & B_k^T \\ 
        B_k & -\kappa^2 \bW_k \end{bmatrix} 
            \begin{bmatrix} \brho_k  \\ 
                \bu_k 
            \end{bmatrix} 
            = \begin{bmatrix} 0  \\ -g\,\bb_k
            \end{bmatrix} ,
        \end{array} \label{eq:linsys_spde}
 \end{equation}
where $M_k$ be the mass matrix associated with inner product $(\brho_k, \s_k)$, $\bW _k$ with the inner product $(u_k, v_k)$ which is diagonal, $B_k$ with the bilinear form $(\div \s_k, u_k)$, and $\bb_k$ is hierarchically generated spatial white noise (generated via Algorithm \ref{alg:hspde}). For a scalable, parallelizable implementation, we have several solvers we may consider, one of which -- hybridization AMG approach from \cite{Lee17,dobrev2019algebraic} -- is amenable to large-scale applications because the mass matrices need only be computed one time (on each level), and then may be reapplied to different realizations of $\bb_k$.

To define the Gaussian densities $\mu_k$ at level $k$ we proceed as in Section \ref{sec:fe_grf}. Let us formally introduce the negative scaled Schur Complement of \eqref{eq:linsys_spde} defined by
$$ A_k :=  \frac{\kappa^2}{g} W_k + \frac{1}{g} B_k M_k^{-1} B_k^T. $$
Then solutions of \eqref{eq:linsys_spde} are Gaussian random vectors with distribution $\mu_k \sim \mathcal{N}(0, A_k^{-1} W_k A_k^{-1})$ and corresponding probability density
$$ \mu_k (u_k) \propto \exp{( - \bu_k^T A_k W_k^{-1} A_k \bu_k) } = \exp{( - \bb_k^T  W_k^{-1}  \bb_k) }.$$
We also define the conditionally Gaussian density $u_k | u_{k+1}$ based on our hierarchical decomposition of white noise in Algorithm \ref{alg:hspde}. 
Sampling from the prior distribution $\mu_k$ and from the conditional distribution $u_k | u_{k+1}$ are summarized in Algorithm \ref{alg:unconditional_sample} and Algorithm \ref{alg:twolevel samples}. 
These algorithms will be used to define the proposal distributions within the multilevel MCMC algorithm in Section \ref{sec:num2}.
\begin{algorithm}
\KwIn{Current level $k$ (with $0 \leq k \leq L$)}
Sample $\bxi_k \sim \mathcal{N}(\boldsymbol{0}, I_k)$\\
Define $ \bb_k = {\bW}_k^{1/2} \boldsymbol{\xi}_k$\\
Compute $\bu_k$ by solving \eqref{eq:linsys_spde}\\
\KwOut{$\bu_k$}  
\protect\caption{Generate a sample $u_k$ from the prior distribution $\mu_k$ at level $k$. \label{alg:unconditional_sample}}
\end{algorithm}

\begin{algorithm}
\KwIn{Current level $k$ (with $0 \leq k < L$), the coarse level sample $\bu_{k+1} = A_{k+1}^{-1}\bb_{k+1}$}
Sample $\boldsymbol{\xi}_k \sim  \mathcal{N}(\boldsymbol{0}, I_k)$\\
Define $ \bb_k =  \bPi^T_k \bb_{k+1} + (I - \bPi^T_k P^T_k) \bW_k^{1/2}\bxi_k$\\
Compute $\bu_k$ by solving \eqref{eq:linsys_spde}\\
\KwOut{$\bu_k$}
\protect\caption{Generate a sample $u_k$ from the conditional distribution $u_k| u_{k+1}$ \label{alg:twolevel samples}}
\end{algorithm}

\subsection{Random Field Realizations Using Hierarchical Components}
\label{secSect:RandomFieldHierarchical}
To visualize the hierarchical components of a fine level solution $u_0$, we consider the Egg model domain \cite{jansen2014egg} using three levels of refinement with 18.5K, 148K, and 1.18M elements for levels $\ell = 2, 1, $ and $0$, respectively. Here we skip over the model details, as these will be addressed in the following section, and focus on the new hierarchical sampler.

\begin{figure}
    \centering
        \includegraphics[trim = 10mm 50mm 30mm 0mm,clip,width=.32\textwidth]{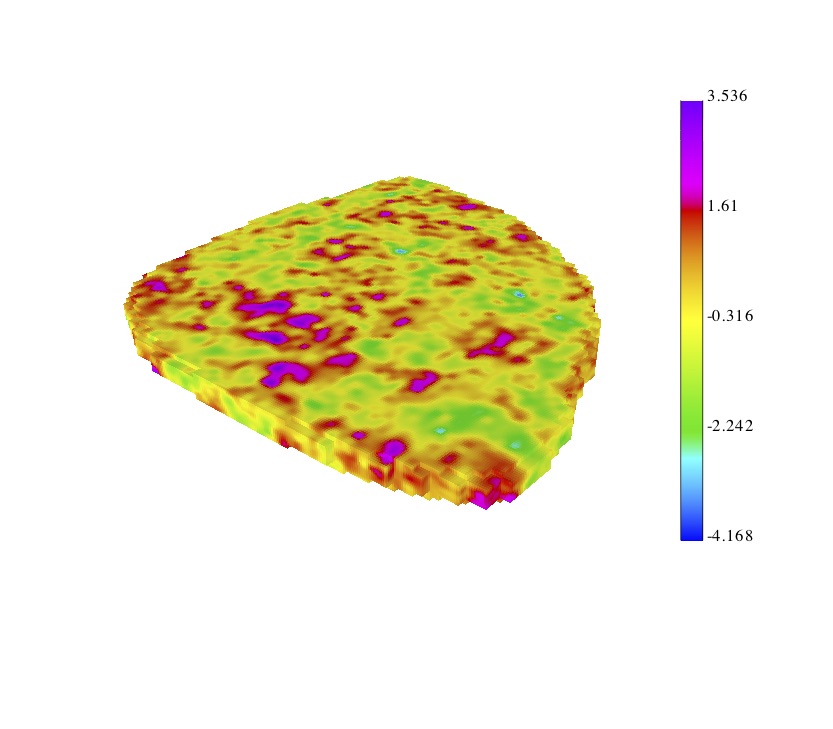}
         \includegraphics[trim = 10mm 50mm 30mm  0mm,clip,width=.32\textwidth]{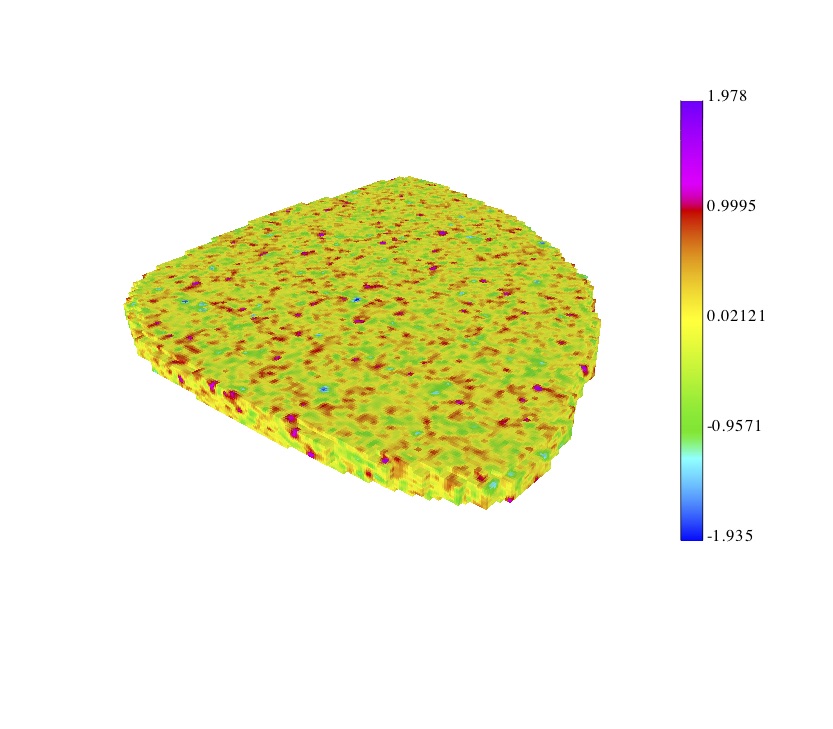}
          \includegraphics[trim = 13mm 50mm 30mm 20mm,clip,width=.32\textwidth]{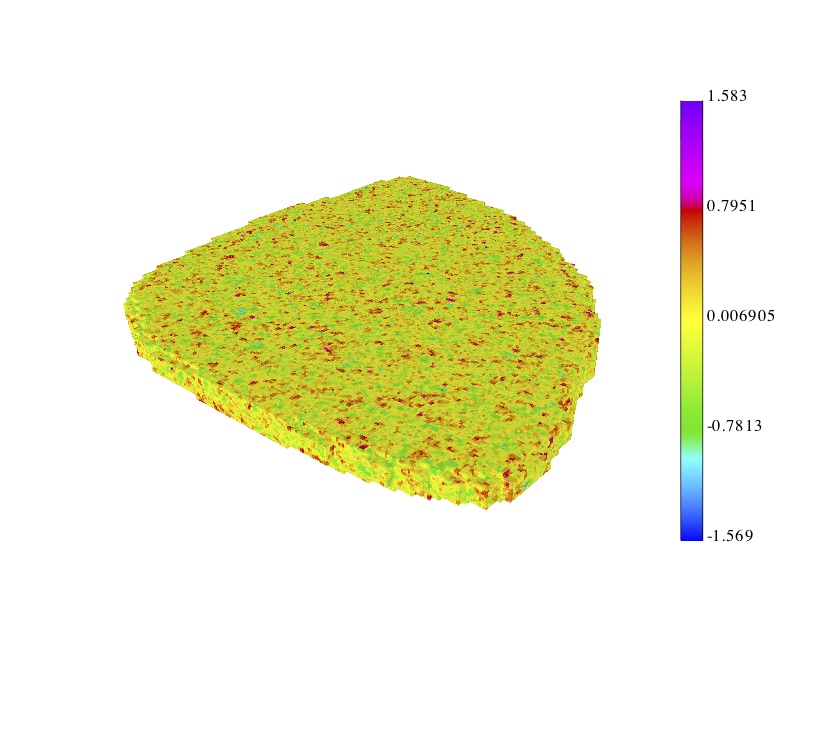}\\
 (a) $u_0^{C2}$  \hspace*{120pt}(b) $u_0^{C1}$ \hspace*{120pt}(c) $u_0^{C0}$ \\
         \includegraphics[trim = 10mm 50mm 30mm 0mm,clip,width=.32\textwidth]{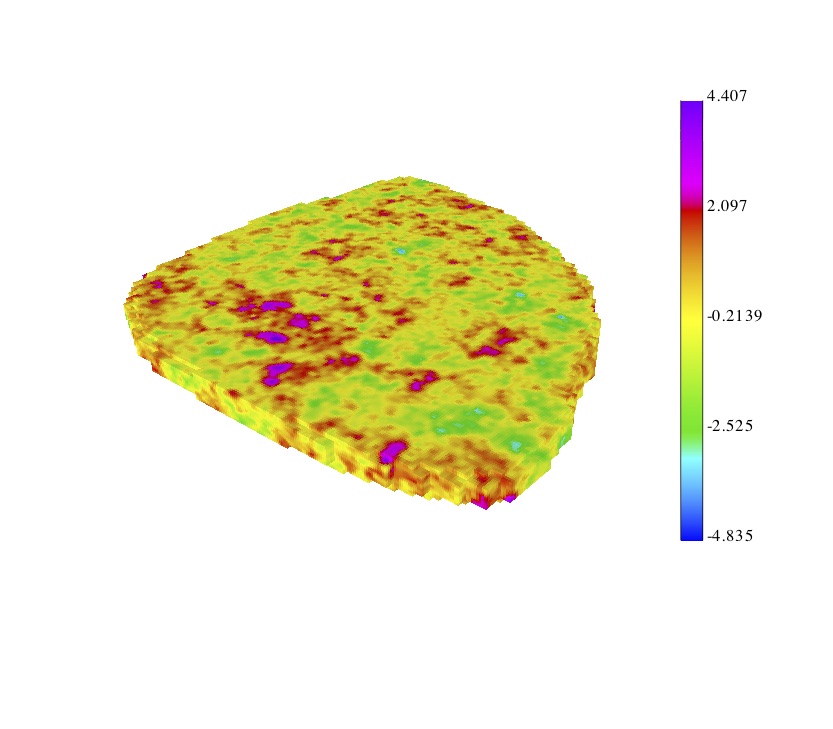}
          \includegraphics[trim = 13mm 50mm 30mm 20mm,clip,width=.32\textwidth]{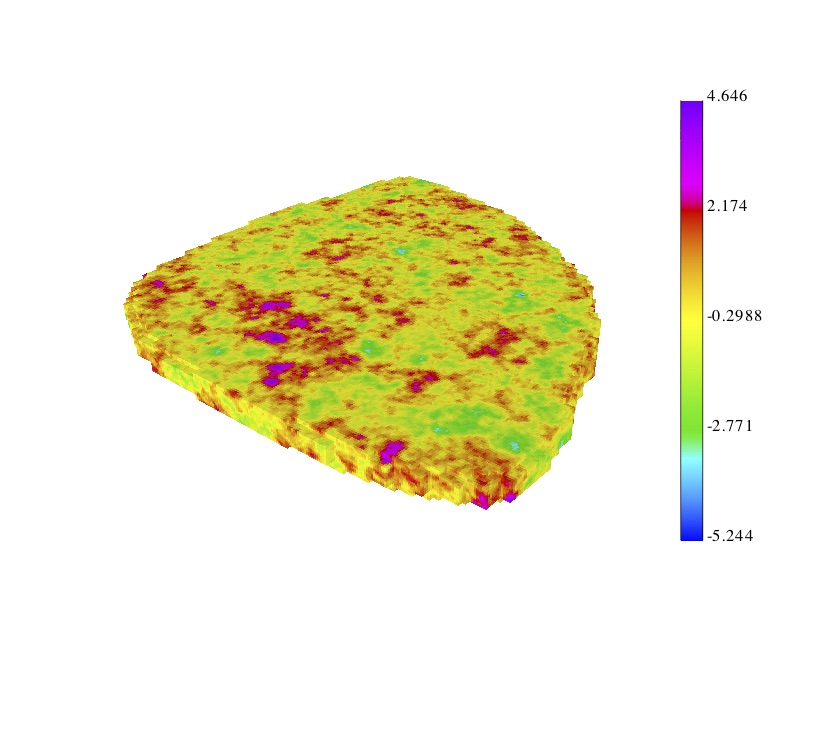}\\
      \hspace*{30pt}     (d) $u_0^{C2}+u_0^{C1}$\hspace*{130pt}(e) $u_0 = u_0^{C2}+u_0^{C1}+u_0^{C0}$\\
   \caption{Various components of a realization $u_0$ on three levels as defined in (\ref{eq:u_comp}), and generated from Algorithm \ref{alg:unconditional_sample} and Algorithm \ref{alg:twolevel samples}. \hrev{Visualizations are rendered with GLVis \cite{glvis-tool}.}}
   \label{fig:egg_grf}
\end{figure}

Using Algorithm \ref{alg:hspde} with $k=0$ and $L=2$, we generate the three components of the righthand side given as $\bb_0 = \Pi_0^T\Pi_1^T W_2^{1/2}\bxi_2 +\Pi_0^T (I - \Pi_1^T P_1^T) W_1^{1/2}\bxi_1 +(I - \Pi_0^T P_0^T) W_0^{1/2}\bxi_0$. For visualization purposes, we separate these three components of $\bb_0$ and solve with each independently. That is, we seek solutions $\bu_0^{C\ell}$ via
\begin{equation}\label{eq:u_comp}
\begin{array}{lcl}
A_0 \bu_0^{C2} &=& \Pi_0^T\Pi_1^T W_2^{1/2}\bxi_2,\\
A_0 \bu_0^{C1} &=& \Pi_0^T (I - \Pi_1^T P_1^T) W_1^{1/2}\bxi_1, \\
A_0 \bu_0^{C0}& =& (I - \Pi_0^T P_0^T) W_0^{1/2}\bxi_0.
\end{array}
\end{equation}
Note that the fine level realization is simply $\bu_0 := \bu_0^{C2}+\bu_0^{C1}+\bu_0^{C0}$.
Figure~\ref{fig:egg_grf}~(a)-(c) displays the solutions $\bu_0^{C2}, \bu_0^{C1}$, and $\bu_0^{C0}$. These results showcase the novelty of this hierarchical approach -- that is, the finite element white noise decomposition enables a realization $u_0$ to be decomposed into independent components across multiple levels. Moreover, this hierarchical approach induces a separation of scales among the terms $u_0^{C\ell}$, similar to that induced by the hierarchical KL-based sampling in \cite{Dodwell15}. On the coarse levels, the terms $u_0^{C\ell}$ capture the smooth components of $u_0$, while, on finer levels, the terms $u_0^{C\ell}$ only contain the highly oscillatory components of $u_0$. This property plays a fundamental role in accelerating the mixing and reducing the variance of the multilevel MCMC in Section \ref{sec:num2}.  This is clearly illustrated by considering the sums of the components  $u_0^{C\ell}$ shown in Figure~\ref{fig:egg_grf}~(d)-(e). In particular, Figure ~\ref{fig:egg_grf}~(d) displays $u_0^{C2}+u_0^{C1}$, and Figure ~\ref{fig:egg_grf}~(e) displays the complete fine level realization $u_0$.    

\section{Numerical Results: Multilevel Hierarchical Sample Generation}\label{sec:num1}
\begin{figure}
    \centering
        \includegraphics[width=.35\textwidth]{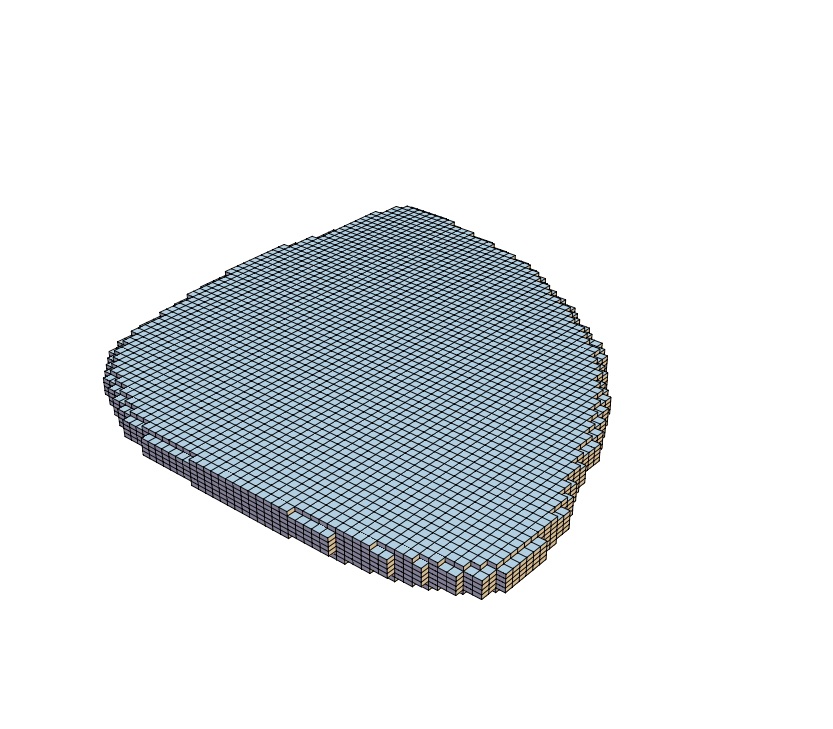}
         \includegraphics[width=.35\textwidth]{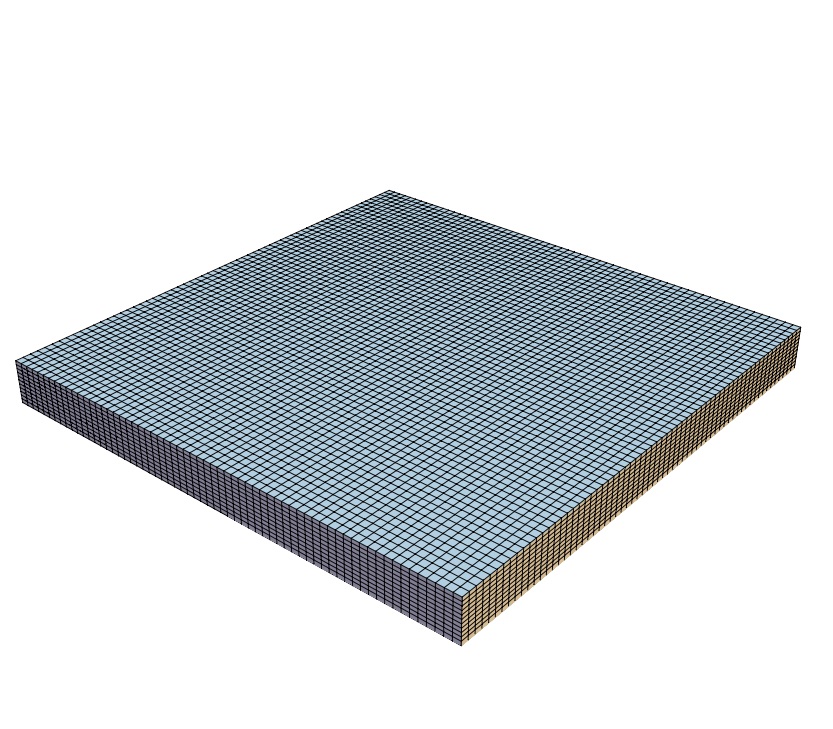}\\
 (a)\hspace*{200pt}(b)
   \caption{(a) Original Egg model mesh containing $18.5$K elements. (b) Enlarged mesh, in which the Egg model mesh is embedded, extends two elements in each direction beyond the Egg model mesh bounding box, and contains $45$K elements. Both meshes displayed correspond to the coarsest level. \label{fig:egg_mesh}}\
  \end{figure}

In this section, we test the hierarchical sampler \hrev{scaling performance} on the `Egg model' \cite{jansen2014egg}, as it contains a large, irregular domain. 
The Egg domain is contained by a $480\text{ m} \times 480\text{ m}\times 28\text{ m}$ bounding box. We note that, as we are employing the approach of \cite{Osborn17b}, we require performing mesh embedding, that is, the Egg model domain is embedded within a $512\text{ m} \times 512\text{ m}\times 44\text{ m}$ domain. This mitigates variance inflation along the boundary due to Neumann boundary conditions (see \cite[Section 2.3]{Lindgren11} and \cite{Osborn17,Osborn17b} for additional discussion). Figure \ref{fig:egg_mesh} displays both the original Egg model mesh and enlarged mesh (in which it is embedded) for the coarsest level, both with hexahedral elements of size $8\text{ m} \times 8\text{ m}\times 4\text{ m}$. Finer mesh resolutions are formed by uniformly refining by a factor of two in each direction.  

We consider three levels $\ell =0,1,2$, with degrees of freedom (corresponding to the number of unknowns in the mixed PDE system as in (\ref{eq:linsys_spde})) given in Table~\ref{tab:egg_weak} with NP $=36$ total MPI processes; then, for a fixed problem size per processor, we increase the number of processes to NP $=288$ and then NP $=2304$. 
\begin{table}  
\caption{Number of global degrees of freedom (DOFs) associated with each level, for each set of process numbers. The DOFs here are associated with the number of unknowns in the mixed PDE system as in (\ref{eq:linsys_spde}).}\label{tab:egg_weak}    
\begin{center}\begin{tabular}{lllllll}
\toprule
  NP &  DOFs $\ell=0$ & DOFs $\ell=1$ & DOFs $\ell=2$  \\
\midrule
  36&  4.8063e+06&  6.0788e+05& 7.7758e+04  \\              
  288&  3.8223e+07& 4.8063e+06& 6.079e+05 \\       
  2304&  3.0488e+08&   3.8223e+07& 4.8063e+06&  \\
\bottomrule
\end{tabular}\end{center}  
\end{table} 
Gaussian random field realizations were generated following our new hierarchical PDE sampling approach; that is, for levels $\ell =0,1,2$, level $\ell$ hierarchical white noise was sampled according to Algorithm \ref{alg:hspde}, and realizations of $u_\ell$ were formed by solving the linear system in (\ref{eq:linsys_spde}) on the Egg domain. 
Numerical simulations were performed using tools developed in ParELAG \cite{parelag}, a parallel C++ library for performing numerical upscaling of finite element discretizations and AMG techniques, and ParELAGMC \cite{parelagmc}, a parallel element agglomeration MLMC library. These libraries use MFEM \cite{mfem-library} to generate the fine grid finite element discretization and HYPRE \cite{hypre} to handle massively parallel linear algebra. \hrev{In particular, we employ hybridization AMG~\cite{Lee17,dobrev2019algebraic}, where the rescaled linear system is solved with conjugate gradient (CG) preconditioned by HYPRE's BoomerAMG}. Note, all timing results were executed on the Quartz cluster at Lawrence Livermore National Laboratory, consisting of 2,688 nodes where each node has two 18-core Intel Xeon E5-2695 processors. For the weak scaling results, we use 36 MPI processes per node. 

\hrev{Table \ref{tab:egg_time} provides the average wall time for the hybridization AMG solver, where the number of preconditioned CG (PCG) iterations -- to reduce the $l^2$ norm of the residual by a factor of $10^6$ -- %
are provided in parentheses. These timing results indicate favorable scaling on the finest level; however, parallel efficiency does degrade on the coarser levels, which is to be expected, as we are limited to our solver and AMG's performance on coarse levels (see~\cite{gahvari2011modeling}). Nonetheless, these results show mesh independence of our approach for all levels in the hierarchy, as the iteration count is stable with increased problem size.}%
\begin{table} 
\caption{Average wall time (seconds) to solve (\ref{eq:linsys_spde}) using hybridization AMG (averaged over $100$ realizations). The average number of PCG iterations are provided in parentheses.}\label{tab:egg_time}
\begin{center}\begin{tabular}{lllll}
\toprule                                                
 Level $\ell$ & Local DOFs  & NP=36   & NP =288 & NP =2304  \\                
 \midrule%
  0     & $135,438$     & 2.42 (11)   & 2.88 (12)  & 3.15 (13)    \\            
  1     &  $17,368$    & 0.187 (10) &  0.258 (11)   &   0.316  (12) \\       
 2      &  $2,280$    & 0.0152  (9) & 0.0297 (11)   & 0.0635 (11)     \\
 \bottomrule                                                                       \end{tabular}\end{center}  
\end{table}  

\begin{figure}
    \centering
        \includegraphics[trim = 0mm 0mm 0mm 0mm,clip,width=.44\textwidth]{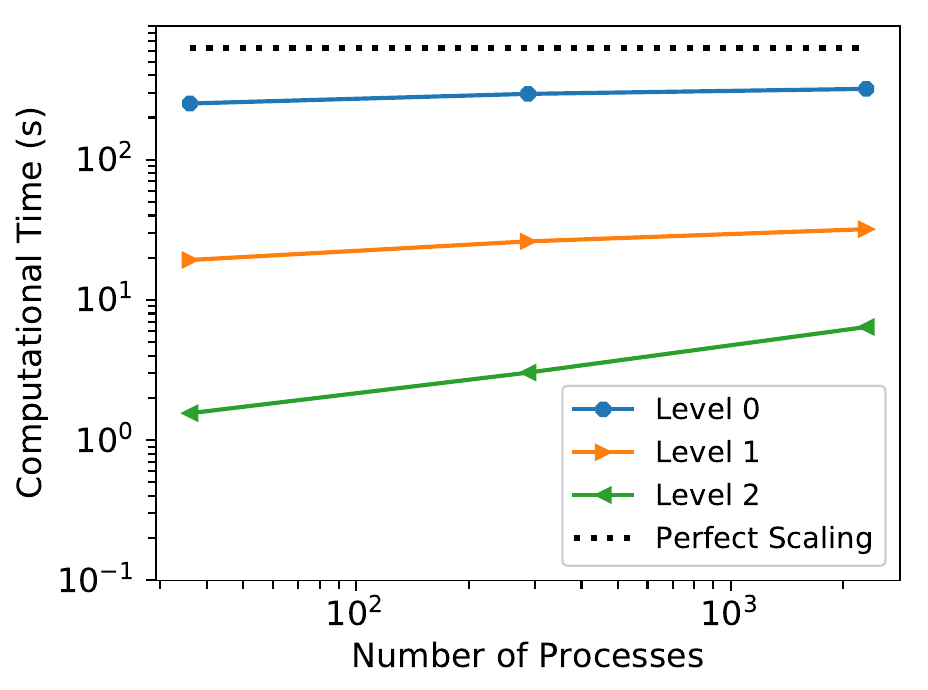}
         \includegraphics[trim = 0mm 0mm 0mm 0mm,clip,width=.44\textwidth]{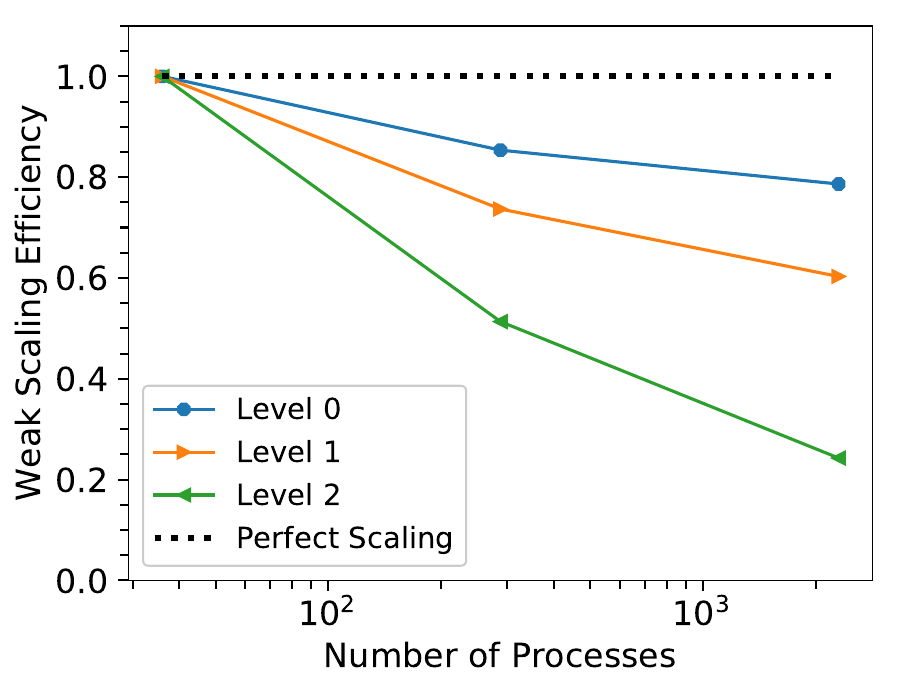}\\
 (a)\hspace*{200pt}(b)
  \caption{(a) Weak scaling for the three different levels; here levels refers to a fixed number of elements per processor. (b) Associated efficiency across the three levels. \label{fig:egg_ws}}
\end{figure}

\hrev{In addition, Figure \ref{fig:egg_ws} displays these weak scaling results (for $100$ simulations), as well as the efficiency decay with increased problem size.}
%
\hrev{While we can claim algorithmic scalability, the degraded scaling on the coarser levels indicates a drawback with our solver implementation, that is, for coarser discretizations we require fewer processes than the fine discretizations to obtain favorable parallel efficiency on all levels. Because of this, we are unable to get complete scaling across all levels; rather each level will benefit from a different number of processes. A topic for a follow up study will include coarse-grid redistribution in order to improve the coarse level performance, see, e.g.,~\cite{trilinos-website}. Nonetheless, it is still clear that the coarse level samples are significantly faster to generate than the finest level samples, which is essential for multilevel MCMC performance.
}

\rsout{In Figure \ref{fig:egg_ws} (a) the computational time for the finest level, $\ell=0$, remains relatively flat with the increase in NP, while the coarsest level, $\ell=2$, has the worst scaling for the given number of processes. That said, it is still clear that the coarse level samples are faster to generate than the finest level samples. For the three levels, the associated efficiency of this scaling is provided in Figure \ref{fig:egg_ws} (b); with the finest level efficiency decreasing to $85\%$ and $80\%$ with the increase in NP, and the coarsest level decreasing to $50\%$ and $25\%$. 
This is to be expected -- with the increase in the number of processes, the coarser levels require greater communication time, which decreases the \hrev{parallel} efficiency; however, it does indicate a drawback with our HPC implementation, that is, for coarser discretizations, we require fewer processes than the fine discretizations. Because of this, we are unable to get complete scaling across all levels; rather each level will benefit from a different number of processes. The work in~\cite{drzisga2017scheduling} describes this as the {\it scalability window} -- that is, each level has a range of NP values it can run on in order to maintain weak scaling efficiencies greater than a specified threshold -- and investigates applying this approach to multilevel Monte Carlo. While this approach would greatly improve the scaling of this method, it is outside the scope of this paper.%
}

\section{Hierarchical PDE Approach for Multilevel MCMC}\label{sec:num2}

In this section, we apply the proposed hierarchical PDE-based sampling approach to solve a nonlinear Bayesian inference problem. We use the multilevel MCMC framework in \cite{Dodwell15} to explore the posterior distribution of the uncertain parameter and estimate moments (the mean) of a scalar quantity of interest $Q$. 

In particular, we consider the problem of inferring a log-normal permeability field from cell-averaged pressure measurements for a single phase steady state subsurface problem. In what follows, we denote with $u \in \Theta$ the uncertain parameter representing the logarithm of the permeability field, with $(\bq, p) \in \mathbf{R}\times \Theta$ the state variables representing the flow velocity and pressure, and with $p_\text{obs} \in \mathbb{R}^m$ the data representing cell-averaged pressure observations at $m$ given measurement locations. %

\emph{Prior distribution.} We assume a Gaussian prior density on the spatially varying log-permeability coefficient, i.e. $u \sim \mathcal{N}(u_*, \mathcal{C})$, with covariance operator $\mathcal{C}$ and mean value $u_* = 0$. To ensure that the inference problem is well-posed in infinite dimensions, we use a squared inverse elliptic operator in \eqref{eq:cov_op} to define the prior covariance operator, see for example \cite{Stuart2010, Flath11, bui2013computational, petra2014computational}.
Samples from the prior distribution can then be drawn by solving \eqref{eq:spde_weak1} as we explained in Section \ref{sec:spde}, and their probability can be computed using \eqref{eq:prior}.

\emph{Forward map.} Let $y = \mathcal{F}(u)$ denote the forward map from the uncertain field $u \in  \Theta$ to the observable $y \in \mathbb{R}^m$.
The map $\mathcal{F} = \mathcal{B}\circ \mathcal{H}$ is the composition of a forward PDE solve $\mathcal{H}$ that computes the pressure field $p$ for a given realization of log-permeability $u$ and a linear operator $\mathcal{B}$ that evaluates the pressure $p$ on local cells.
For a single phase porous media flow, with $k = \exp(u)$, we define $p = \mathcal{H}(u)$ as the solution to the mixed formulation of the Darcy's equations
    \begin{equation}
       \begin{array}{ll}
            (k^{-1}\bq , \s)- (\div \s, p)= (\f , \s) & \forall \s \in \mathbf{R} \\
            (\div \bq, v) = 0 & \forall v \in \Theta,
        \end{array}\label{eq:darcy_weak}
    \end{equation}
with Dirichlet boundary condition $p = p_D$ on $\Gamma _D$, enforced by the right-hand side $\f$, and Neumann boundary condition $\bq \cdot \bn = 0$ on $\Gamma_N$, where $\Gamma_D$ and $\Gamma _N$ are non-overlapping partitions of $ \partial D$. %

\emph{Likelihood function.} We assume that the measured data $p_\text{obs} \in \mathbb{R}^m$ are corrupted by additive Gaussian noise $\eta$ with zero mean and covariance $\Gamma_\eta = \sigma_\eta^2 I_m$, where $I_m$ is the identity matrix in $\mathbb{R}^m$. That is,
\begin{equation}\label{eq:yobs}
p_\text{obs}= \mathcal{F}(u) + \eta, \quad \eta \sim \mathcal{N}(0, \Gamma_{\eta}).
\end{equation}
From the noise model in \eqref{eq:yobs}, we have that the conditional probability of $p_\text{obs}$ given $u$ is also Gaussian with mean $\mathcal{F}(u)$ and covariance  $\Gamma_{\eta}$, that is
$$ p_\text{obs} | u \sim \mathcal{N}(\mathcal{F}(u), \Gamma_{\eta} ).$$
The likelihood function  $\pi^{\rm like}(p_\text{obs} | u)$ then reads
\begin{equation}\label{eq:like}
\pi ^\text{like}( p_\text{obs}| u  ) \propto \exp \left( -\frac{1}{2}\Vert p_\text{obs} - \mathcal{F}(u) \Vert^2_{\Gamma_\eta^{-1}} \right),
\end{equation}
where $\| \cdot \|^2_{\Gamma_{\eta}^{-1}}$ denotes the $\Gamma_{\eta}^{-1}$-weighted $l^2$ norm in $\mathbb{R}^m$.

\emph{Posterior distribution.} By applying Bayes' theorem, the posterior density $\nu$ in the infinite dimensional case  is given by
\begin{equation}\label{eq:post}
\nu (u | p_\text{obs}) \propto \pi ^\text{like}( p_\text{obs}| u  ) d\mu(u),
\end{equation}
where $\pi ^\text{like}( p_\text{obs}| u  )$ is the likelihood function in \eqref{eq:like} and $d\mu(u)$ is the prior density in \eqref{eq:prior}.
We note that, although the prior distribution and likelihood functions are both Gaussian, the posterior distribution $\nu(u | p_\text{obs})$  is not Gaussian because of the nonlinearity introduced by the forward map $\mathcal{F}$. Thus, there is no closed form solution to the Bayesian inference problem and therefore we will use MCMC sampling to explore the posterior distribution.

\emph{Quantity of interest.} Finally, let us introduce the scalar quantity of interest $Q = Q(u)$ representing the flux across the outflow boundary $\Gamma_\text{out}$, which is defined as
\begin{equation}\label{eq:Q}
Q = \frac{1}{|\Gamma_\text{out}|}\int_{\Gamma_\text{out}} \bq (\cdot, \omega)\cdot \bn \ dS,
\end{equation}
where $\bn$ represents the outward unit vector normal to $\Gamma_\text{out} \subset \partial D$.

Our goal is to estimate the posterior mean of $Q$, defined as

\begin{equation}\label{eq:postQ}
\E_\nu [Q] = \int_\Omega Q(u) \, d \nu(u | p_\text{obs}) ,
\end{equation}
by sampling the posterior distribution \eqref{eq:post} using multilevel MCMC. As a reference to the reader, notation introduced and frequently used in this section is provided in Table \ref{tab:notation2}.
\begin{table}[!]  \footnotesize
\caption{Bayesian inference notation.}\label{tab:notation2}
\begin{center}\begin{tabular}{ll}
\toprule 
 Infinite-Dimensional Variable   & Description\\                              
\midrule
$ d\mu (u)$ & Prior density of $u$ \\
$ p_\text{obs} $& Observed local pressure data in $D$ \\
 $\pi^{like}(p_\text{obs}|u)$ & Likelihood function \\
 $\nu (u|p_\text{obs})$ & Posterior density \\
$\E_\nu [\cdot] $ & Mean with respect to posterior density \\
 $Q=Q(u)$ & Quantity of interest \\
\midrule
Finite-Dimensional Variable   & \\  
\midrule
$ d\mu _\ell (u_\ell)$ & Prior density of $u_\ell$ on level $\ell$ \\
 $\pi_\ell^{like}(p_\text{obs}|u_\ell)$ & Likelihood function \\
 $\nu_\ell (u_\ell|p_\text{obs})$ & Posterior density on level $\ell$ \\
 $\alpha_\ell^{SL} $ & Single-level acceptance probability on level $\ell$ \\
  $\alpha_\ell^{ML} $ & Multilevel acceptance probability on level $\ell$ \\
$\E_{\nu_\ell} [\cdot] $ & Mean with respect to level $\ell$ posterior density\\
$\V_{\nu_\ell} [\cdot] $ & Variance with respect to level $\ell$ posterior density\\
 $Q_\ell=Q(u_\ell)$ & Quantity of interest on level $\ell$ \\
\midrule
   \bottomrule                                                                         
    \end{tabular}\end{center}  
\end{table}  

\subsection{Markov Chain Monte Carlo}
For the single-level approach, the log-permeability $u$, associated likelihood $\pi ^\text{like}( p_\text{obs}| u  )$, and QoI $Q$ (see (\ref{eq:spde_weak1}), (\ref{eq:like}), (\ref{eq:Q}), and, respectively) are approximated numerically by solving the mixed PDEs in (\ref{eq:spde_weak1}) and (\ref{eq:darcy_weak}) utilizing a finite element approach on triangulation $\TT_\slk$ (for the finest level $\slk$). We denote these discrete approximations as $u_\slk $, $\pi^{\text{like}}_\slk$, and $Q_\slk = Q(u_\slk)$. %

MCMC, and in particular, Metropolis-Hastings, is a modified Monte Carlo approach, where samples $Q_\slk$ are generated from the target (posterior) distribution via a Markov chain. Then the posterior QoI expectation $\E_{\nu_k}[Q]$ may be approximated as
\begin{equation}\label{eq:mcmcQ}
\hat{Q}_\slk^{MCMC} =\frac{1}{N} \sum\limits_{i=n+1}^{n+N} Q_\slk^{(it)},
\end{equation}
where $Q_\slk^{(it)}= Q(u_\slk^{(it)})$, $n$ is the number of samples discarded as burn-in, and $t$ is the subsampling rate to obtain independent samples (see Appendix \ref{sec:iact}). 
To generate subsequent samples within a chain, a new $u_\slk^{\text{prop}}$ is sampled from the proposal distribution and subjected to Metropolis-Hastings acceptance/rejection criterion. In this work, we utilize a preconditioned Crank-Nicolson stepping scheme with step size $\beta>0$, where samples $\bm \psi_{\slk}$ from $\mu_{\slk}$ are drawn using Algorithm \ref{alg:unconditional_sample}. Thanks to the prior-invariance of the preconditioned Crank-Nicolson proposal, the sample $u_\slk^{\text{prop}}$ is then accepted with probability $\alpha_\slk^{SL}$ defined in (\ref{eq:alpha_sl}) \cite{cotter2013mcmc}. The procedure is summarized in Algorithm \ref{alg:q_sl}; additional details can be found in \cite{Dodwell15}.
\begin{algorithm}[!]
\begin{itemize}
\item Given  $u_\slk^{(i-1)}$, propose $u_{\slk}^{\text{prop}}$ using preconditioned Crank-Nicolson:
\begin{equation}\label{eq:pcn1}
u_{\slk}^{\text{prop}} := \sqrt{1-\beta^2} u _{\slk}^{(i-1)} + \beta \bm \psi_{\slk},
\end{equation}
where $\bm \psi_{\slk} \sim \mu_{\slk}$ is computed using Algorithm \ref{alg:unconditional_sample}.
 \item Accept $u_\slk^{(i)} = u_{\slk}^{\text{prop}}$ with probability 
 \begin{equation}\label{eq:alpha_sl}
 \alpha_\slk^{SL} (u_\slk^\text{prop}|u_\slk^{(i-1)}) = \min \left\{ 1, \frac{\pi ^\text{like}_\slk( p_\text{obs}| u_\slk^\text{prop}  ) }{\pi ^\text{like}_\slk( p_\text{obs}| u_\slk^{(i-1)}  )}  \right\}
 \end{equation}
 \item Return $u_\slk^{(i)}$ and $Q_\slk^{(i)} = Q_\slk(u_\slk^{(i)})$.
\end{itemize}
\protect\caption{Single level Metropolis-Hastings MCMC Algorithm with preconditioned Crank-Nicolson proposal to generate a posterior sample $u_\slk^{(i)} | u_\slk^{(i-1)}$\label{alg:q_sl}}
\end{algorithm}

The cost of performing MCMC depends on how quickly the chain mixes as well as the variance of $\hat{Q}_\slk$. The first---the mixing of the chain---is controlled by the autocorrelation of samples within the chain. As adjacent samples in the chain are correlated (and not independent), the integrated autocorrelation time $\tau_Q$ of the chain will determine how many steps (and thus forward simulations) are required to get to the next independent sample. 
The second---the variance of $\hat{Q}_\slk$---is controlled by the number of independent samples used in the estimator, i.e., $N$. The accuracy is similar to Monte Carlo in that the required number of (independent) samples to achieve a desired mean squared error depends on the variance of $Q_\slk$ as well as the bias introduced by numerically approximating $Q$. If we require $N$ independent simulations, with an integrated autocorrelation time (rounded up to an integer value) of $t$ (see Appendix \ref{sec:iact}), then we require at least $t N$ simulations. Thus acceleration approaches should seek to reduce $t$ and $N$ by increasing the mixing of the chain and reducing the variance of the estimator, respectively.

\subsection{Multilevel Markov Chain Monte Carlo}
To accelerate MCMC, we consider the multilevel framework in \cite{Dodwell15}, which utilizes chains at coarser spatial discretization levels to perform the majority of likelihood functions evaluations. Similar to previous sections, let us denote the log-normal permeability field, the Darcy pressure and flux, and the QoI at dicretization level $\ell$ with the symbols $u_\ell$, $(p_\ell, \bq_\ell)$, $Q_\ell := Q_\ell(u_\ell)$, respectively, for $k=0\leq\ell\leq L$.  Then the posterior mean of $Q_0$ is equivalently written as
\begin{equation}\label{eq:postQML}
\E_{\nu_0} [Q_0] = \E_{\nu_L} [Q_L]  + \sum\limits_{\ell=0}^{L-1} \left(\E_{\nu_\ell} [Q_\ell]  - \E_{\nu_{\ell+1}} [Q_{\ell+1}]\right) ,
\end{equation}
where $\nu_{\ell}$ is the discrete posterior measure on level $\ell$. 
Following \cite{Dodwell15}, for each level $\ell$, we define a multilevel estimator $\hat{Y}_\ell^{N_\ell}$ of the difference $\E_{\nu_\ell} [Q_\ell]  - \E_{\nu_{\ell+1}} [Q_{\ell+1}]$ and write
\begin{equation}\label{eq:postY}
\hat{Y}_\ell^{N_\ell} =\frac{1}{N_\ell} \sum\limits_{i=n_\ell+1}^{n_\ell+N_\ell} Y_\ell^{(it_{\ell})} = \frac{1}{N_\ell} \sum\limits_{i=n_\ell+1}^{n_\ell+N_\ell} \left(Q_\ell^{(it_{\ell})}-Q_{\ell+1}^{(it_{\ell}t_{\ell+1})}\right).
\end{equation}
Above $n_\ell$ corresponds to the burn-in on level $\ell$, $N_\ell$ is the effective sample size on level $\ell$ (defined later in (\ref{eq:nleff})), $t_\ell$ and $t_{\ell+1}$ are the estimated integrated autocorrelation times of the chains at levels $\ell$ and $\ell+1$, respectively (see Appendix \ref{sec:iact} for details).
The key aspect of the multilevel MCMC is how to couple Markov chains at different levels so that: \emph{i}) the variance of $Y_\ell$ is much smaller than that of $Q_\ell$, \emph{ii}) information from coarser level chains are used to accelerate mixing of finer level chains (higher acceptace rate, smaller integrated autocorrelation time). 
In this section, we focus on how to generalize the multilevel MCMC algorithm \cite[Algorithm 3]{Dodwell15} to replace the KL decomposition-based sampling with our scalable multilevel PDE samplers described in Section \ref{sec:hspde_implementation}.

As motivated in Remark \ref{rem:2level}, in what follows, we describe a two-level chain to evaluate the difference estimator $\hat{Y}_\ell$ at a generic level $0 \leq \ell < L$.  Given a coarse sample $u_{\ell+1}^{(j-1)t_{\ell+1}}$, we advance the coarse chain at level $\ell+1$ by $t_{\ell+1}$ steps using single-level Metropolis-Hastings as in Algorithm \ref{alg:q_sl}. This results in a coarse sample $u_{\ell+1}^{(jt_{\ell+1})}$ that is independent of $u_{\ell+1}^{(j-1)t_{\ell+1}}$.

To propose $u_{\ell}$ on the finer level $\ell$, we use the two-level preconditioned Crank-Nicolson in \eqref{eq:pcn2}, where $\bm \psi_{\ell}$ is sampled from the conditional distribution $\bm \psi_{\ell} | u_{\ell+1}^{(j-1)t_{\ell+1}}$ using Algorithm  \ref{alg:twolevel samples}. Note that the independence of $u_{\ell+1}^{(jt_{\ell+1})}$ from $u_{\ell+1}^{(j-1)t_{\ell+1}}$ guarantees that also $\bm \psi_{\ell}$ is independent of $u_{\ell}^{(j-1)}$. That means that the two-level preconditioned Crank-Nicolson proposal in \eqref{eq:pcn2} satisfies the assumptions of \cite[Lemma 3.1]{Dodwell15}, and therefore the multilevel acceptance probability $\alpha_{\ell}^{ML} (u_\ell^\text{prop}|u_\ell^{(i-1)})$ in \eqref{eq:alpha_ml} satisfies the detailed balance condition. Algorithm \ref{alg:q_ml} summarizes the generation of the paired fine and coarse level samples $u_\ell^{(j)} $ and $u_{\ell+1}^{(jt_{\ell+1})}$.\\
Note that, if $u_\ell^{\text{prop}} $ is accepted at step $j$, then $u_\ell^{(j)} $ and $u_{\ell+1}^{(jt_{\ell+1})}$ are correlated. Specifically, both $u_\ell^{(j)} $ and $u_{\ell+1}^{(jt_{\ell+1})}$ are generated from the same coarse level white noise functional $\bb_{\ell+1}$, and thus the difference $u_\ell^{(j)} - P_\ell u_{\ell+1}^{(jt_{\ell+1})}$ is small. This is observed in Section \ref{secSect:RandomFieldHierarchical}, where Figure~\ref{fig:egg_grf}~(b)-(c) display these differences, defined as solutions $u_0^{C\ell}$ (see (\ref{eq:u_comp})).
This means that one should expect $Y_\ell^{(j)}$ to be small when step $j$ is accepted, which is a necessary condition to achieve multilevel acceleration. The numerical results presented next demonstrate that our algorithm is indeed able to achieve multilevel acceleration of the chain mixing and variance reduction.

\begin{algorithm}[!]
\underline{Part I. Advance the coarse chain at level $\ell+1$ by $t_{\ell+1}$ steps:}\\
\begin{itemize}
\item Given $u_{\ell+1}^{(j-1)t_{\ell+1}}$, apply Algorithm \ref{alg:q_sl} on level $\ell+1$ for $t_{\ell+1}$ steps
\item Store $u_{\ell+1}^{(jt_{\ell+1})}$ and $Q(u_{\ell+1}^{(jt_{\ell+1})})$
\end{itemize}
\underline{Part II. Advance the fine chain at level $\ell$ by one step}:\\
\begin{itemize}
\item Given $u_{\ell}^{(j-1)}$ and $u_{\ell+1}^{(jt_{\ell+1})}$, propose $u_{\ell}^{\text{prop}}$:
 \begin{equation}\label{eq:pcn2}
 u_{\ell}^{\text{prop}} := \sqrt{1-\beta^2}u_{\ell}^{(j-1)} + \beta \bm \psi_\ell,
 \end{equation}
where $\bm \psi_\ell | u_{\ell+1}^{(jt_{\ell+1})}$ is sampled  using Algorithm \ref{alg:twolevel samples}.
\item Accept $u_{\ell}^{(j)} = u_{\ell}^{\text{prop}}$ with probability
\begin{equation}\label{eq:alpha_ml}
\alpha_{\ell}^{ML} (u_\ell^\text{prop}|u_\ell^{(i-1)}) = \min \left\{ 1, \frac{\pi ^\text{like}_\ell( p_\text{obs}| u_\ell^\text{prop}  )\pi ^\text{like}_{\ell+1}( p_\text{obs}| u_{\ell+1}^{(j-1)t_{\ell+1}} ) }{\pi ^\text{like}_\ell( p_\text{obs}| u_\ell^{(j-1)}  )\pi ^\text{like}_{\ell+1}( p_\text{obs}| u_{\ell+1}^{jt_{\ell+1}}  )}  \right\}.
\end{equation}
\item Return  $u_{\ell}^{(j)}$, $u_{\ell+1}^{(jt_{\ell+1})}$, $Y_\ell^{(j)} = Q(u_\ell^{(j)}) -Q(u_{\ell+1}^{(jt_{\ell+1})})$ 
\end{itemize}
\protect\caption{Two-level Metropolis-Hastings MCMC Algorithm to generate paired samples  $u_{\ell+1}^{(jt_{\ell+1})}$ and $u_{\ell}^{(j)} | u_{\ell}^{(j-1)} $ \label{alg:q_ml}}
\end{algorithm}

\begin{remark}\label{rem:2level}
Another small difference with respect to the work in \cite{Dodwell15} is that each $\hat{Y}_\ell^{N_\ell}$ estimate uses only two levels. Specifically, Algorithm \ref{alg:q_ml} uses a single auxiliary chain on the coarser level $\ell+1$ to estimate $\hat{Y}_\ell$, while \cite{Dodwell15} runs auxiliary chains on all coarser levels. Our decision to do this is based on algorithmic simplicity and scalability. That said, utilizing all coarser levels is a detail that may be considered in future work.
\end{remark}

\subsection{Four-Level Markov Chain Monte Carlo Results}

To demonstrate our method \hrev{is well-suited for multilevel MCMC algorithms,} we test a \hrev{four}-level MCMC with the hierarchical stochastic PDE solver by utilizing Algorithm \ref{alg:q_ml}. The computational domain is a unit cube discretized using tetrahedral elements, \hrev{with about $1.57$M elements on the finest level. Each coarser level is formed by uniformly coarsening by a factor of $8$, until obtaining $3,072$ elements on the coarsest level.} For the prior, we consider Gaussian random fields with correlation length $\lambda = 0.3$ and marginal variance of $\sigma^2 =0.5$. 
\rsout{Table \ref{tab:cube_time} provides information regarding the number of elements in each level, as well as the wall time needed to generate $50$ realizations of $u_k$ for various numbers of processors.  
%
\begin{table}[!]  
\caption{Time in seconds to generate 50 Gaussian random field realizations using (\ref{eq:linsys_spde}). Runs were completed on LLNL Quartz.}\label{tab:cube_time}
\begin{center}\begin{tabular}{llllllll}
\toprule                                                
 Level $\ell$ &  Number Elements & NP=2 & NP =4 & NP =8 & NP =16 & NP =32 & NP = 64 \\                              \midrule
  0&  196,608 & 308& 161   & 75.2      & 42.9      & 23.2  &  12.0\\                 
  1&  24,576  & 27.9 &  10.9    &   5.87    &  3.88    & 2.94&  2.55 \\       
 2&  3,072  &  2.18 & 1.49    & 1.11    &  1.03      &   1.07   & 1.18  \\            \bottomrule                                                                          \end{tabular}\end{center}  
\end{table}  
While we may consider a level dependent number of processors based on these scaling results (which is a possible future direction), in the following computations we fix the total number of processors to $8$.} For the observational pressure data, we synthetically generate a realization $p_\text{obs}\in\R^{25}$ with $\sigma_\eta^2 =0.005$; this is done via \eqref{eq:yobs} on a reference mesh with approximately $12.5$M elements. 

\hrev{We run five independent two-level chains (as in Algorithm \ref{alg:q_ml}) using step size $\beta^2 = 0.3$ to estimate $\hat{Y}_\ell$ at levels $\ell=0,1,2$. 
For each of these chains, the first $1,000$ samples of $Y_\ell$ are discarded as burn-in; then the following $1,000$ independent samples---properly subsampled according to integrated autocorrelation time estimates---are used in our statistical approximations. Similarly we estimate $\hat{Q}_3$ with five independent single-level chains (as in Algorithm \ref{alg:q_sl}); however on this coarsest level, we use a longer burn-in of $3,000$ samples.}

\begin{figure}[!]
    \centering
        \includegraphics[width=.3\textwidth]{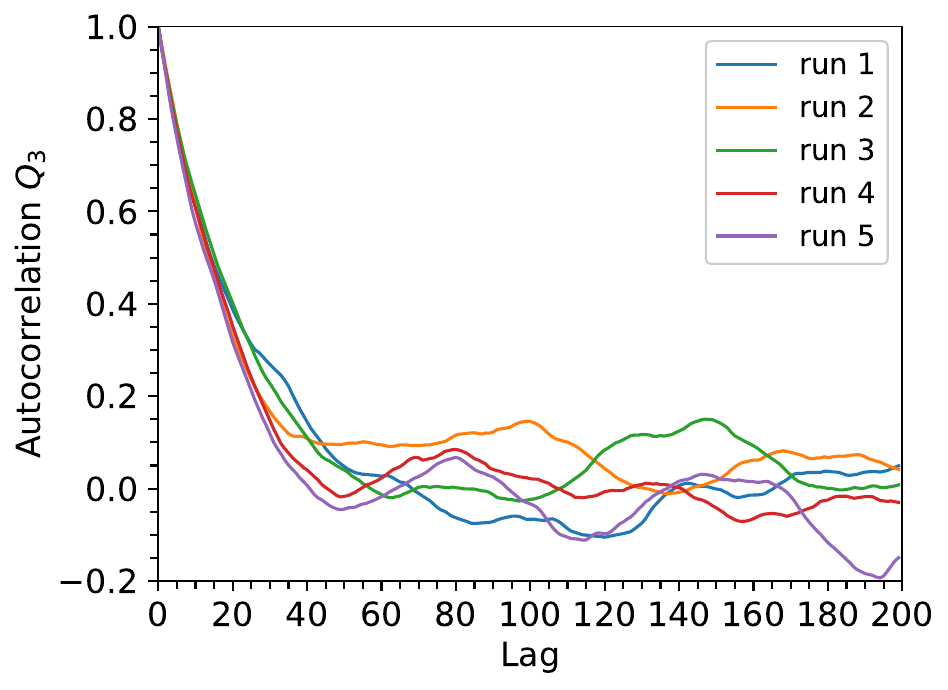}
        \includegraphics[width=.3\textwidth]{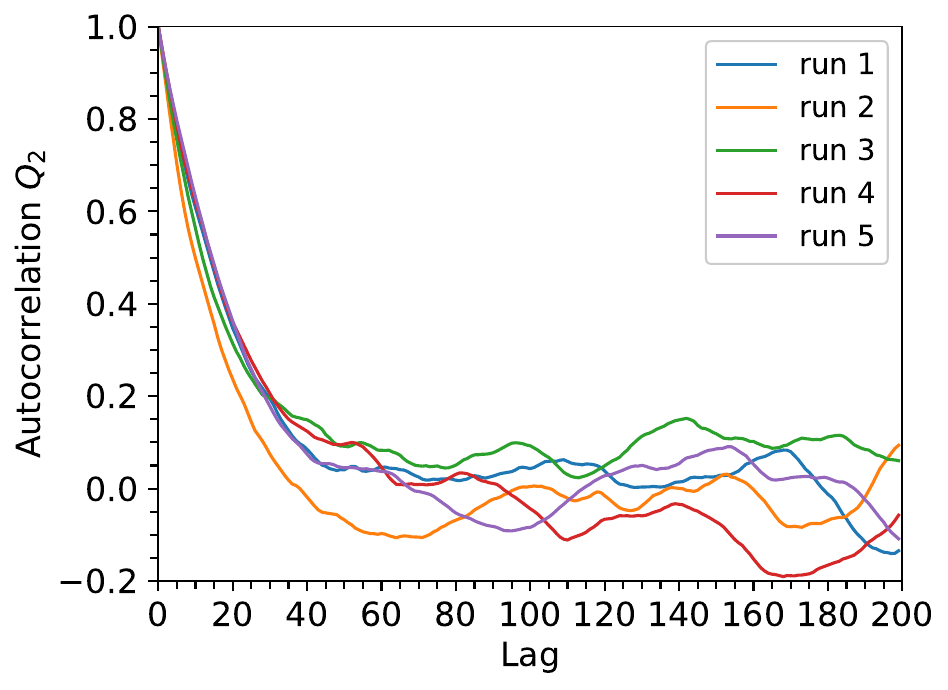}
         \includegraphics[width=.3\textwidth]{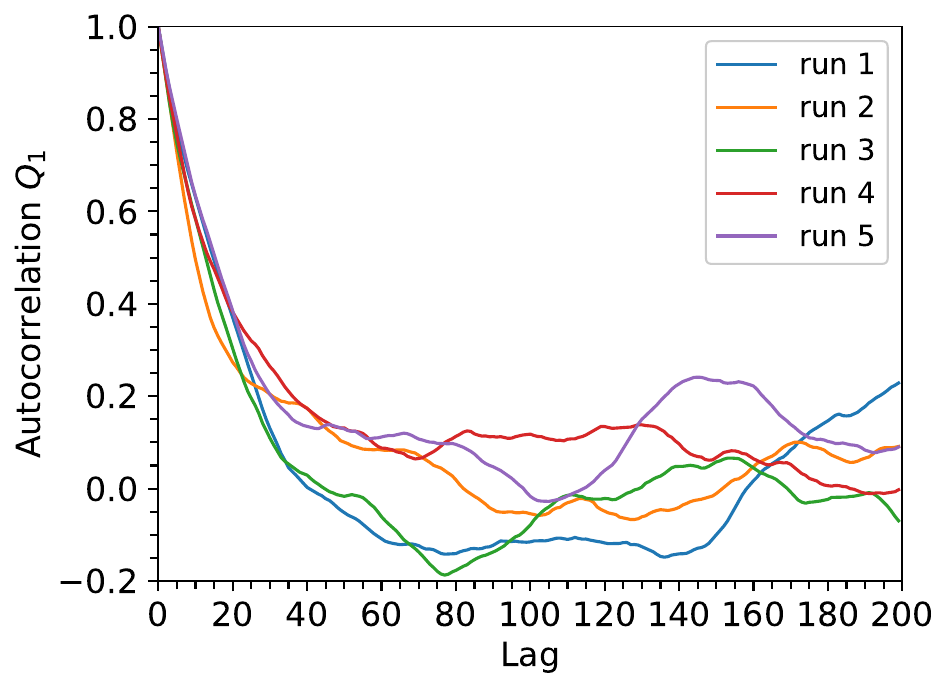}\\
        \includegraphics[width=.3\textwidth]{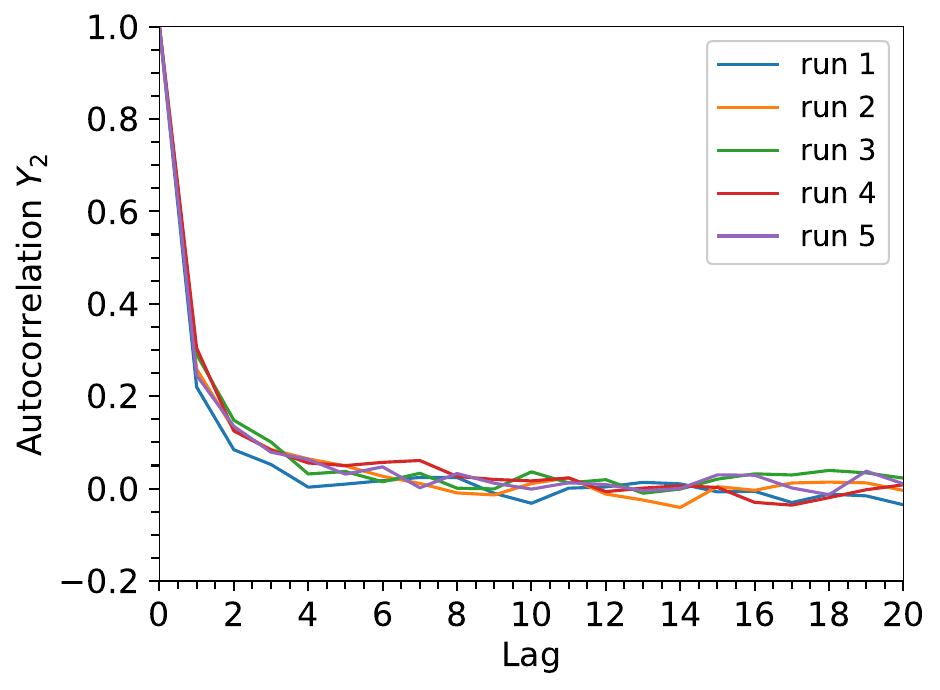}
        \includegraphics[width=.3\textwidth]{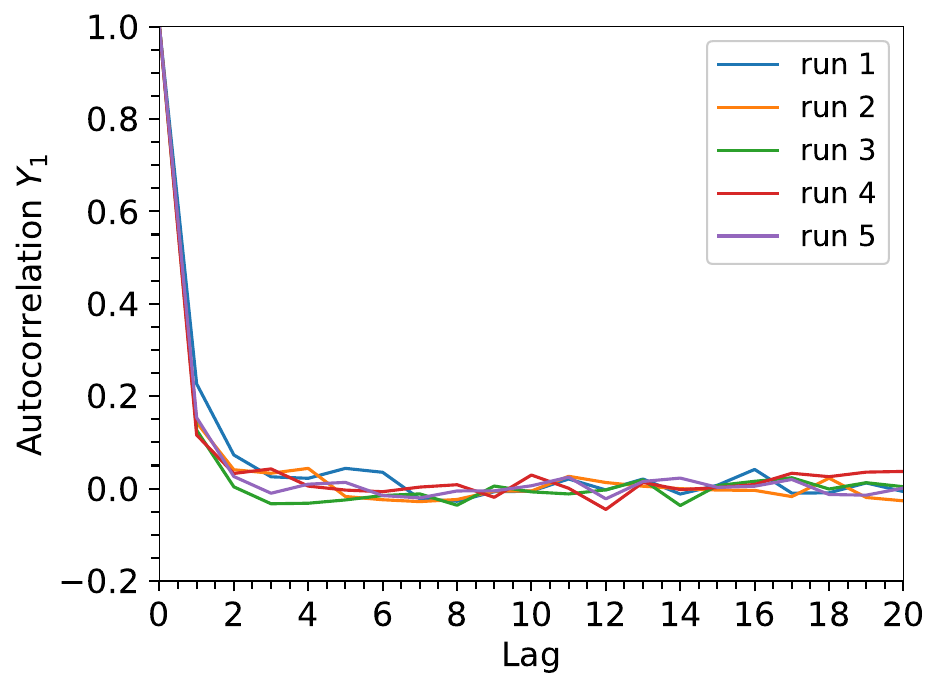}
         \includegraphics[width=.3\textwidth]{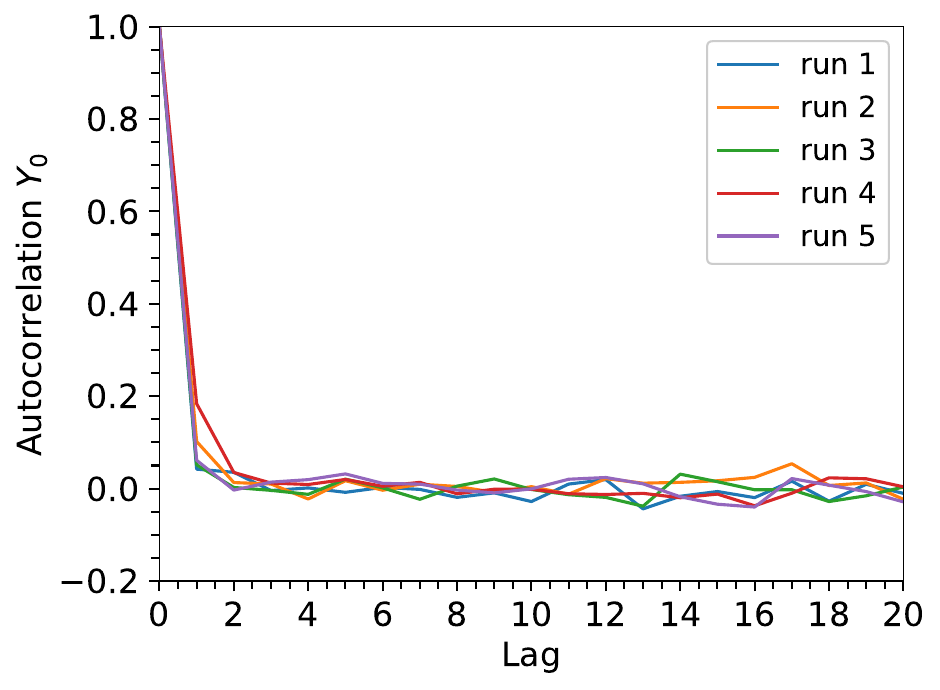}\\
   \caption{\hrev{Autocorrelation estimates for increasing lag time for 5 two-level chains used to estimate $Y_\ell=Q_\ell - Q_{\ell+1}$ at levels $\ell=0,1,2$. Top row displays the autocorrelation for the coarse $Q_{\ell+1}$ samples, while the bottom row shows the autocorrelation for corresponding correction chain $Y_\ell$.} \label{fig:cube_autocorr}}
\end{figure}

\begin{figure}[!]
\centering
        \includegraphics[trim = 0mm 0mm 0mm 0mm,clip,width=.4\textwidth]{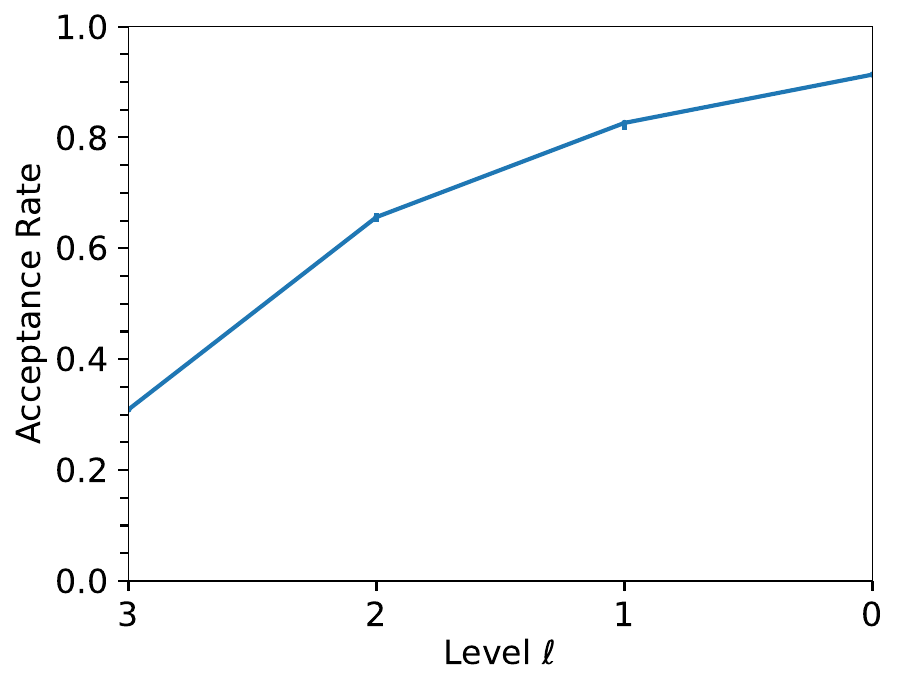}
         \includegraphics[trim = 0mm 0mm 0mm 0mm,clip,width=.41\textwidth]{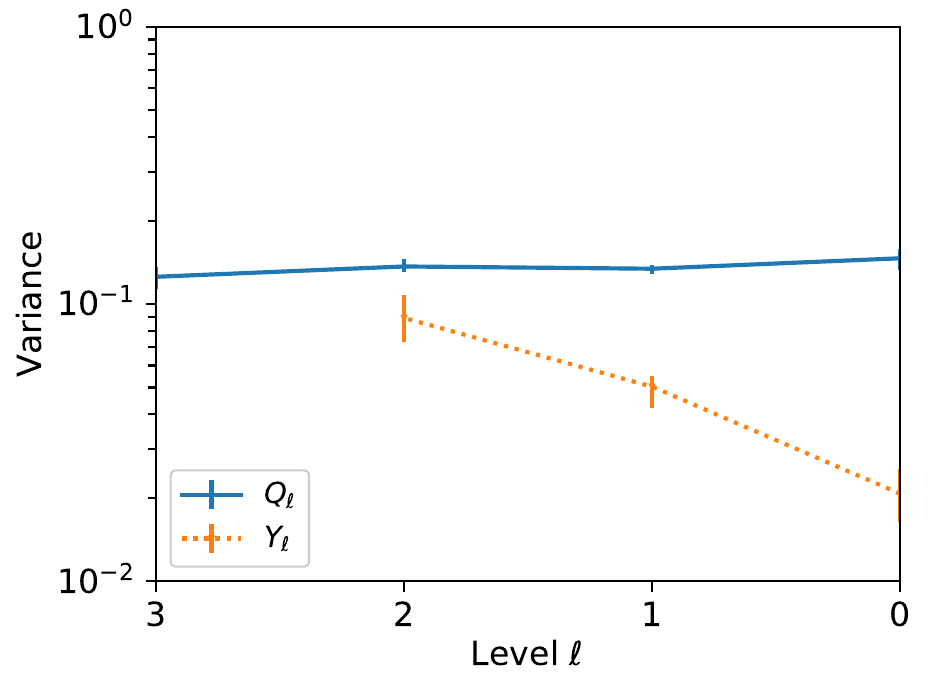}
          \hspace*{30pt}(a)\hspace*{180pt}(b)
        \caption{(a) Average Acceptance rate for the two-level MCMC Algorithm \ref{alg:q_ml} as a function of the level $\ell$. (b) Average variance estimates for $Q_\ell$ and $Y_\ell$. \hrev{Averages are taken over 5 independent runs, with error bars indicating minimum and maximum values.}} 
        \label{fig:cube_accept}
\end{figure}

\hrev{From the five resulting chains of $Y_\ell = Q_{\ell}-Q_{\ell+1}$ (for each $\ell=0,1,2$), we estimate the autocorrelation as a function of lag time. In Figure \ref{fig:cube_autocorr}, the top row displays these estimates for $Q_{\ell+1}$, which indicate the mixing of the coarse chain in the estimation of $Y_\ell$. The bottom row of Figure \ref{fig:cube_autocorr} displays these estimates for $Y_{\ell}$, where the faster decay in autocorrelation is the result of the coarse level being subsampled. We further note that the autocorrelation time of $Y_\ell$ decays faster for smaller $\ell$, as we expect the acceptance rate to increase with mesh refinement as shown in Figure \ref{fig:cube_accept}~(a). }

Figure~\ref{fig:cube_accept}~(a) displays the average acceptance rate from the five chains, on the four different levels. \hrev{The small error bars indicate the range of acceptance rate values from the five chains.} This increase in the acceptance rate \hrev{with the refinement of levels} is similar to that reported in \cite{Dodwell15}, numerically demonstrating that indeed our method of using multilevel stochastic PDE samplers is a computationally efficient alternative to KL-decomposition based sampling for multilevel MCMC. 
\rsout{Next, we consider the estimated variance of the $Y_\ell$ for levels $\ell=0,1,2$.} 
Figure~\ref{fig:cube_accept}~(b) provides average variance estimates for $Q_\ell$ on each level, as well as average variance estimates for the correction terms $Y_\ell = Q_\ell-Q_{\ell+1}$. \hrev{Error bars indicate the range of variance values from the five chains.} The decay in the $Y_\ell$ variance estimate indicates that fewer samples are required on the finer levels (relative to the coarsest) to obtain a target mean square error. This result is similar to that of \cite{Dodwell15}; however, our decay is not quite as rapid. This is likely due to the difference in problem setup ($d=3$, $\lambda = 0.3$), as well as the fact that we're doing inference in a higher-dimensional space. More specifically, we have one DOF per element with $196$K elements on the finest level, while the work in \cite{Dodwell15} has a KL expansion truncated at $150$ DOFs on the finest level. \hrev{Furthermore, we note this decay is dependent on the multilevel acceptance rates for each level. In Algorithm \ref{alg:q_ml}, upon rejection of a fine sample, the realization $Y_\ell^{(j)} = Q(u_\ell^{(j)}) -Q(u_{\ell+1}^{(jt_{\ell+1})})$ is calculated from unrelated realizations of $u_{\ell}$ and $u_{\ell+1}$. This feature results in a variance that decays slower than in typical multilevel Monte Carlo (see, e.g., \cite{Giles08, Giles13, Cliffe11, Teckentrup13}). }

\begin{table}[!]  
\caption{\hrev{Multilevel MCMC average estimates from 5 chains. Values of the effective cost and sample size are calculated via (\ref{eq:cleff}) and (\ref{eq:nleff}), with $\varepsilon=0.01$}. }\label{tab:cube_stats}
\begin{center}\begin{tabular}{llllllllllr}
\toprule
Level $\ell$&  $\mathcal{C}_\ell$&  $t_{\ell+1}$&  $t_\ell$&  $\mathcal{C}_\ell^{\text{eff}}$&  $\E[Q_\ell]$&  $\E[Y_\ell]$&  $\E[\vert Y_\ell\vert]$&  $\V[Q_\ell]$&  $\V[Y_\ell]$&  Predicted $ N_\ell^{\text{eff}}$ \\ 
\midrule
0&  494.08&  44&  4&  12902.4&  1.22&  0.0054&  0.0366&  0.15&  0.0207&  773 \\ 
1&  62.08&  36&  3&  1033.0&  1.20&  0.0046&  0.0808&  0.13&  0.0503&  4265 \\ 
2&  7.84&  40&  5&  239.2&  1.21&  0.0169&  0.1509&  0.14&  0.0891&  11788 \\ 
3&  1.00&  {$-$}&  45&  45.0&  1.17&  1.1709&  1.1709&  0.13&  0.1254&  32245 \\ 
\bottomrule
\end{tabular}\end{center}
\end{table}  

\hrev{In Table~\ref{tab:cube_stats}, we provide the averaged statistical estimates derived from these five chains (on each level). In particular, we provide the estimated integrated autocorrelation times for each level, as well as the estimated mean and variance values for $Q_\ell$ and $Y_\ell$. We note that, while the $\V[Y_\ell]$ and $\E[\vert Y_\ell\vert]$ approximations decay with mesh refinement, the $\E[Y_\ell]$ approximations are not monotonic. However, this result does not conflict with the expected results, as the variance indicates error in our estimates. Using these multilevel estimates of $E[Y_\ell]$ and $\E[Q_3]$ we predict a posterior mean of $1.198$, with a variance (based on the number of effective samples defined in (\ref{eq:nleff})) of $5\cdot 10^{-5}$. With an equivalent cost, we expect the single-level estimator to have a variance of $1.74\cdot 10^{-4}$.}

\hrev{Finally we compare the predicted cost to run this approach with that of single-level MCMC, based on the calculations in Table \ref{tab:cube_stats}. As the scalability of solvers was investigated in \cite{Lee17,dobrev2019algebraic,Osborn17b}, which show cost per iteration is proportional to the DOFs, and we demonstrate the number of iterations stays stable with mesh refinement, we define the cost per simulation on level $\ell$, denoted $\mathcal{C}_\ell$, to be the number of global DOFs in our linear system, normalized with respect to the coarsest level. Subsequently,} the \urev{optimal} effective cost per independent sample of $Y_\ell$ is defined as
\begin{equation}\label{eq:cleff}
\mathcal{C}_\ell^{\text{eff}}:= t_\ell (\mathcal{C}_\ell + t_{\ell+1} \mathcal{C}_{\ell+1})
\end{equation}
for $\ell=0,\hdots, L-1$, and $\mathcal{C}_\ell^{\text{eff}} = t_\ell \mathcal{C}_\ell$ for $\ell=L$. 
\hrev{Then, from \cite{Dodwell15},} the effective sample size (for a mean square error tolerance of $\varepsilon^2$) on each level is calculated via 
\begin{equation}\label{eq:nleff}
N_\ell^{\text{eff}} = \frac{2}{\varepsilon^2}\left(\sum\limits_{k=0}^L \sqrt{ \V _{\nu_k,\nu_{k+1}}[Y_k] \mathcal{C}_k^{\text{eff}}} \right)\sqrt{ \frac{\V _{\nu_\ell,\nu_{\ell+1}} [Y_\ell] }{\mathcal{C}_\ell^{\text{eff}}}},
\end{equation}
where $\V_{\nu_\ell,\nu_{\ell+1}} [Y_\ell]$ is the variance of $Y_\ell$ with respect to the joint distribution of $u_\ell$ and $u_{\ell+1}$. \hrev{For single-level MCMC, the number of effective samples for the target mean square error $\varepsilon^2$ is $N_{sl}^\text{eff}= 2\V_{\nu_0}[Q_0]/\varepsilon^2$, and the effective cost is $\mathcal{C}_{sl}^\text{eff}= t_L \mathcal{C}_0$. Using these results, we estimate that performing our four-level MCMC is about $3.5$ times faster than the single-level approach. In comparison the hierarchical four-level approach of \cite{Dodwell15} is about $5$ times faster than the single-level approach. Aside from different statistical estimates, a key difference of these approaches (that impacts the cost) lies in the number of levels used to estimate each $Y_\ell$ chain. As noted before, the approach of \cite{Dodwell15} utilizes all coarser levels in a hierarchical manner to estimate these differences. Although this component of the algorithm was not implemented in our numerical results, it is a feature that we would like to include in future work.}

\section{Conclusion}
\hrev{In this work we develop a novel, (algorithmically) scalable, hierarchical PDE-based approach to generate Gaussian random field realizations} that is well-suited for \hrev{multilevel MCMC on} large-scale three-dimensional problems.\rsout{Building off of approaches from \cite{Osborn17,Osborn17b}, we extend the work into a hierarchical multilevel framework by performing a hierarchical decomposition of white noise across multiple levels of discretization. From this hierarchical decomposition, we may generate corresponding realizations of discrete Gaussian random fields, such that fine scale random fields are generated in a hierarchical fashion from coarse scale random fields.} 
%
The novelty and advantages of our two-level preconditioned Crank-Nicolson proposal in Algorithm \ref{alg:q_ml} lies in the use of a scalable, memory efficient stochastic PDE sampler in lieu of a computationally and memory expensive KL-decomposition\rsout{. Our method shares the same structure and benefits of the original KL-based approach} in \cite{Dodwell15}. Similarly to the proposals in \cite{Dodwell15}, our proposals are linear transformations of independent Gaussian vectors defined on the coarse and fine grids: the coarser-level random variables define the smooth components of the random field $u_\ell$, while the finer-level random variables control the high frequency components of $u_\ell$. However,  our method uses sparse finite element interpolation operators and scalable fast PDE linear solvers to define such linear transformation, while the one in \cite{Dodwell15} uses dense matrices whose columns represents the dominant eigenvectors of the covariance method.

As our numerical result showed, Algorithm \ref{alg:q_ml} offers comparable multilevel acceleration to that presented in \cite{Dodwell15}. 
First, the great majority of likelihood evaluations are done the coarse levels of the hierarchy, where evaluating the forward model is inexpensive. Second, the acceptance rate improves as the mesh is refined thus reducing the variance of the estimator $\hat{Y}_\ell$ at finer levels.  Third, the auxiliary coarse level chain allows for drastically reducing the integrated autocorrelation time $t_\ell$ thanks to the use of independent samples from the coarse chain. 
As numerically illustrated in Section \ref{secSect:RandomFieldHierarchical}, our hierarchical sampler induces a multiscale decomposition of the random field $u$, where the finer-level proposal $u_\ell$ shares the same smooth components of the corresponding sample from the posterior distribution $\nu_{\ell+1}$ at the coarser-level $\ell+1$. As we move to finer and finer levels we expect the likelihood function to become insensitive to the difference $u_\ell - Pu_{\ell+1}$, thus drastically increasing the acceptance rate. The increased mixing of the chain is then a direct consequence of the increased acceptance rate and of the independence of the coarse grid samples used in the two-level preconditioned Crank-Nicolson proposal.

\hrev{The next stage of research will include investigating the overall scaling of multilevel MCMC with this new hierarchical sampler. In particular, an important---and necessary---component will be coarse grid redistribution for improved performance on all levels in the sampling hierarchy.}
In addition, possible future directions of this work include performing this multilevel MCMC approach with a derivative enhanced proposal, e.g., local Hessian information, as in \cite{cui2019multilevel}, which combines the multilevel approach of \cite{Dodwell15} with dimension-independent likelihood-informed MCMC samplers of \cite{cui2016dimension} to further accelerate multilevel MCMC.

\rsout{A current drawback of this implementation is that the coarser levels do not scale as well. This is due to the number of processors; in particular, the fine level scales in a different regime than the coarse levels, and thus we see decreased efficiency on the coarser levels, as shown in the numerical results. To improve this issue, one approach would be to utilize the framework of \cite{drzisga2017scheduling}, where the authors distribute the work across a number of processors that is level dependent. }

\rsout{Possible future directions of this work include performing this multilevel MCMC approach with a more informed proposal, e.g., local Hessian information, as in \cite{cui2016dimension, cui2019multilevel}. In particular, \cite{cui2019multilevel} combines the multilevel approach of \cite{Dodwell15} with dimension-independent likelihood-informed MCMC samplers of \cite{cui2016dimension} to further accelerate multilevel MCMC. }

\section*{Acknowledgements}
This document was prepared as an account of work sponsored by an agency of the United States government. Neither the United States government nor Lawrence Livermore National Security, LLC, nor any of their employees makes any warranty, expressed or implied, or assumes any legal liability or responsibility for the accuracy, completeness, or usefulness of any information,
apparatus, product, or process disclosed, or represents that its use would not infringe privately owned rights. Reference herein to any specific commercial product, process, or service by trade
name, trademark, manufacturer, or otherwise does not necessarily constitute or imply its
endorsement, recommendation, or favoring by the United States government or Lawrence
Livermore National Security, LLC. The views and opinions of authors expressed herein do not
necessarily state or reflect those of the United States government or Lawrence Livermore
National Security, LLC, and shall not be used for advertising or product endorsement purposes. %

%
\bibliographystyle{plain}
\bibliography{HF_bib}

\appendix
\section{Integrated Autocorrelation Time}\label{sec:iact}
To obtain independent samples for unbiased estimates of QoI moments from the chain $\{Q_0^{(i)}\}_{i>0}$, we subsample the chain according to its integrated autocorrelation time $\tau_Q$. In this work, we estimate $\tau_Q$ as 
\begin{equation}\label{eq:iact}
\hat{\tau}_Q = 1+ 2\sum\limits_{\tau= 1}^{M}\hat{\rho}_Q(\tau)
\end{equation}
where the normalized autocorrelation function is estimated as 
\begin{equation}\label{eq:act}
\hat{\rho}_Q(\tau)= \frac{1}{N-\tau}\sum\limits_{i = 1}^{N-\tau}\frac{(Q_0^{(i)} - \hat{\mu}_Q)(Q_0^{(i+\tau)} - \hat{\mu}_Q)}{\hat{\sigma}_Q^2},
\end{equation}
with $ \hat{\mu}_Q$ and $\hat{\sigma}_Q^2$ as the estimated mean and variance (respectively) of the data $\{Q_0^{(i)}\}_{i=1}^N$, and $M \ll N$ (see \cite{sokal1997monte} for more information on integrated autocorrelation time). 

In practice we subsample at a rate of $t :=\ceil{\hat{\tau}_Q}$. We denote $t_\ell$ (with $\ell<L$) as the estimate for the multilevel chains $\{Y_\ell^{(i)}\}_{i>0}$.

\end{document}